\documentclass[12pt]{amsart}
\usepackage{amssymb}

\usepackage{graphicx}

\usepackage[]{amsfonts}

\def\underset#1#2{\mathrel{\mathop{\kern0pt #2}\limits_{#1}}}

\def\overset#1#2{\mathrel{\mathop{\kern0pt #2}\limits^{#1}}}

\def\couleur(#1 #2 #3)
	{
\special{ps::[end]}}

\def\sqr#1#2{{\vcenter{\vbox{\hrule height.#2pt
             \hbox{\vrule width.#2pt height#1pt \kern#1pt
             \vrule width.#2pt}
             \hrule height.#2pt}}}}

\def\st{\mathinner{\mkern1mu\raise1pt\hbox{.}				% three points : "such that"
		   \mkern1mu\raise4pt\hbox{.}
		   \mkern1mu\raise1pt\hbox{.}
		 }
         }

\def\bx#1{\setbox1=\hbox{\kern3pt{#1}\kern3pt}				% Make a box. Close it by "}"
 \dimen1=\ht1 \advance\dimen1 by 3pt \dimen2=\dp1 \advance\dimen2 by 3pt
 \setbox1=\hbox{\vrule height\dimen1 depth\dimen2\box1\vrule}%
 \setbox1=\vbox{\hrule\box1\hrule}%
 \advance\dimen1 by .4pt \ht1=\dimen1
 \advance\dimen2 by .4pt \dp1=\dimen2 \box1\relax}

\def\k#1{\kern#1em}
\def\vci{\vrule  width.02em height1.47ex depth-.0ex}				% le 1 en blackboard
\def\11{{\rm\k{.2}\vci\k{-.37}1}}

\newtheorem{Theorem}{Theorem}[section]
\newtheorem{Proposition}[Theorem]{Proposition}
\newtheorem{Remark}[Theorem]{Remark}

\begin{document}
\title{Numerical solution of a parabolic problem arising in finance}
\author{Marie-Noelle Le Roux}
\address{UNIVERSITE BORDEAUX1, Institut de Math{\'e}matiques de Bordeaux, 
UMR 5251,351,Cours de la Lib{\'e}ration, 33405, Talence Cedex}
\email{Marie-Noelle.Leroux@math.u-bordeaux1.fr}
\maketitle
\begin{abstract} {
In this paper, we study a parabolic system of three equations which permits 
to solve an optimal replication problem in incomplete markets. We obtain 
existence and uniqueness of the solution in suitable Sobolev spaces and 
propose a numerical method to compute the optimal strategy.\ \par
}\end{abstract}
\ \par
\section{Introduction}
\setcounter{equation}{0}We study here a parabolic system arising in the 
resolution of an optimal replication problem in incomplete  markets. 
Given a European derivative security with an arbitrary payoff function, 
the optimal replication problem is to find a dynamic portfolio strategy, 
that is self-financing and comes as close as possible to the payoff at 
maturity date $T$. In complet markets, such a dynamic-hedging strategy 
exists: the payoff of a European option can be replicated exactly; it 
is the Black-Scholes model (1973)~\cite{Black}.In ~\cite{Kogan} , Bertsimas, 
Kogan and Lo propose a solution approach for this problem in incomplete 
markets..\ \par
At time $\tau =0$, consider a portfolio of stocks and riskless bonds 
at a cost $V_{0}$ and denote by $\theta (\tau ),\ B(\tau ),\ V(\tau )$ 
the number of shares of the stock held, the value of bonds held and the 
market value of the portfolio at time $\tau $. Hence, $\displaystyle 
V(\tau )=\theta (\tau )P(\tau )+B(\tau ),\ 0\leq \tau \leq T$.\ \par
If we note $\sigma $ the volatility and $F$ the payoff function, the 
value function $J$ is defined by: 
\begin{displaymath} 
J(\tau ,V,P,\sigma )=\underset{\theta (s),\ s\geq \tau }{\min }E(((V(T)-F(P(t),\sigma 
(T)))^{2}/(V(\tau ),P(\tau ),\sigma (\tau ))).\end{displaymath}  The 
replication error $\epsilon (V_{0})$ is $\left({J(0,V,P,\sigma )}\right) ^{1/2}$ 
and it can be minimized  with respect to the initial wealth $V_{0}$ to 
yield the least-cost optimal-replication strategy and the minimum replication 
error $\epsilon ^{*}$ is $\epsilon ^{*}=\underset{V_{0}}{{\mathrm{min}}}\epsilon (V_{0})$\ 
\par
In ~\cite{Kogan}, it has been proved the the value function  $J$ is quadratic 
in $V:$ $J=a(V-b)^{2}+c$ and the coefficients  $a,\ b,\ c$ satisfy the 
following system of partial differential equations: \ \par

\begin{equation} 
\frac{\displaystyle \partial a}{\displaystyle \partial \tau }=-\frac{\displaystyle 
k^{2}\sigma ^{2}}{\displaystyle 2}\frac{\displaystyle \partial ^{2}a}{\displaystyle 
\partial \sigma ^{2}}-g_{1}(\sigma )\frac{\displaystyle \partial a}{\displaystyle 
\partial \sigma }+\rho ^{2}k^{2}\frac{\displaystyle \sigma ^{2}}{\displaystyle 
a}\displaystyle \left({\displaystyle \frac{\displaystyle \partial a}{\displaystyle 
\partial \sigma }}\right) ^{2}+af^{2}(\sigma ),\label{MarchesIncomplets1}
\end{equation} \ \par
\ \par

\begin{displaymath} 
\frac{\partial b}{\partial \tau }=-\frac{k^{2}\sigma ^{2}}{2}\frac{\partial 
^{2}b}{\partial \sigma ^{2}}-\frac{\sigma ^{2}P^{2}}{2}\frac{\partial 
^{2}b}{\partial P^{2}}-\rho k\sigma ^{2}P\frac{\partial ^{2}b}{\partial 
\sigma \partial P}-g_{2}(\sigma )\frac{\partial b}{\partial \sigma }\end{displaymath} 
\ \par

\begin{equation} 
-(1-\rho ^{2})k^{2}\frac{\displaystyle \sigma ^{2}}{\displaystyle a}\frac{\displaystyle 
\partial a}{\displaystyle \partial \sigma }\frac{\displaystyle \partial 
b}{\displaystyle \partial \sigma },\label{MarchesIncomplets2}
\end{equation} \ \par
\ \par

\begin{displaymath} 
\frac{\partial c}{\partial \tau }=-\frac{k^{2}\sigma ^{2}}{2}\frac{\partial 
^{2}c}{\partial \sigma ^{2}}-\frac{\sigma ^{2}P^{2}}{2}\frac{\partial 
^{2}c}{\partial P^{2}}-\rho k\sigma ^{2}P\frac{\partial ^{2}c}{\partial 
\sigma \partial P}-g(\sigma )\frac{\partial c}{\partial \sigma }-\sigma 
f(\sigma )P\frac{\partial c}{\partial P}\end{displaymath} \ \par
\ \par

\begin{equation} 
-(1-\rho ^{2})k^{2}\sigma ^{2}a\displaystyle \left({\displaystyle \frac{\displaystyle 
\partial b}{\displaystyle \partial \sigma }}\right) ^{2},\label{MarchesIncomplets3}
\end{equation} \ \par
\ \par
where $k>0,$ $\rho \in [-1,+1]$, ( $\rho $ is a correlation coefficient)\ 
\par
$g(\sigma )=-\delta \sigma (\sigma -\sigma _{1})$ ( $\delta >0$ and $\sigma _{1}\in ]0,1[$)\ 
\par
\ \par
$\displaystyle g_{1}(\sigma )=g(\sigma )-2\rho k\sigma f(\sigma ),g_{2}(\sigma )=\ g(\sigma 
)-\rho k\sigma f(\sigma )$\ \par
\ \par
and $\displaystyle f(\sigma )=\displaystyle \left\lbrace{\displaystyle \begin{matrix}{\displaystyle 
\frac{\displaystyle \mu }{\displaystyle \sigma _{0}}\ if\ \ \sigma \leq 
\sigma _{0}}\cr {\displaystyle \frac{\displaystyle \mu }{\displaystyle 
\sigma }\ if\ \sigma \geq \sigma _{0}}\cr \end{matrix}}\right. ,\ \mu 
>0,\ (\mu $ is the drift).\ \par
\ \par
\begin{Remark}  The function $f$ has been modified near $0$ in order 
to be bounded and assure the existence of a solution.\ \par
\end{Remark}
\ \par
The conditions at the time expiry $T\ $ are given by: $\displaystyle 
a(T)=1,\ b(T)=F(P,\sigma ),\ c(T)=0$.\ \par
Under the optimal replication strategy $\theta $$^{{\mathrm{*}}}$, the 
minimum replication error as a function of the initial wealth $V_{0}\ $is 
$(J(0))^{\frac{1}{2}}=\left({a(0)(V_{0}-b(0))^{2}+c(0)}\right) ^{\frac{1}{2}}$, 
hence the initial wealth that minimizes the replication error is $V_{0}^{*}=b(0)\ $ 
the minimal replication error over all $V_{0}$ is $\epsilon ^{*}=\sqrt{c(0)}$ 
and the least--cost optimal strategy at $\tau =0$ is $\displaystyle \theta ^{*}(0)=\frac{\displaystyle \partial b}{\displaystyle \partial 
P}(0)+\frac{\displaystyle \rho k}{\displaystyle P}\frac{\displaystyle 
\partial b}{\displaystyle \partial \sigma }(0)$.\ \par
\ \par
\begin{Remark} Exact replication is possible when $k^{2}(1-\rho ^{2})=0$ 
and this corresponds to the following cases:\ \par
{\hskip 1.8em}- Volatility is a deterministic function of time.\ \par
{\hskip 1.8em}- The Brownian motions driving stocks prices and volatility 
are perfectly correlated.\ \par
\end{Remark}
\ \par
\ \par
{\hskip 1.8em}In this paper, we propose a numerical method to compute 
the solution of equations ~(\ref{MarchesIncomplets1}), ~(\ref{MarchesIncomplets2}), 
~(\ref{MarchesIncomplets3}) and then obtain the minimal replication error 
and the least-cost optimal replication strategy. \ \par
{\hskip 1.8em}To obtain a forward problem, we change the sense of time; 
we note $t=T-\tau $. In order to avoid the function $a$ at the denominator, 
we make the change of unknown $u_{1}={\mathrm{ln}}(a)$. We also replace 
$\sigma $ by $x$, $P$ by $y$, $b$ by $u_{2}$ and $c$ by $u_{3}$.\ \par
The preceding system becomes:  \ \par

\begin{equation} 
\frac{\displaystyle \partial u_{1}}{\displaystyle \partial t}-\frac{\displaystyle 
k^{2}x^{2}}{\displaystyle 2}\frac{\displaystyle \partial ^{2}u_{1}}{\displaystyle 
\partial x^{2}}-g_{1}(x)\frac{\displaystyle \partial u_{1}}{\displaystyle 
\partial x}+k^{2}(\rho ^{2}-\frac{\displaystyle 1}{\displaystyle 2})x^{2}\displaystyle 
\left({\displaystyle \frac{\displaystyle \partial u_{1}}{\displaystyle 
\partial x}}\right) ^{2}+f^{2}(x)=0,\label{MarchesIncomplets0}
\end{equation} \ \par
\ \par

\begin{displaymath} 
\frac{\partial u_{2}}{\partial t}-\frac{k^{2}x^{2}}{2}\frac{\partial ^{2}u_{2}}{\partial 
x^{2}}-\frac{x^{2}y^{2}}{2}\frac{\partial ^{2}u_{2}}{\partial y^{2}}-\rho 
kx^{2}y\frac{\partial ^{2}u_{2}}{\partial x\partial y}-g_{2}(x)\frac{\partial 
u_{2}}{\partial x}\end{displaymath} \ \par

\begin{equation} 
-(1-\rho ^{2})k^{2}x^{2}\frac{\displaystyle \partial u_{1}}{\displaystyle 
\partial x}\frac{\displaystyle \partial u_{2}}{\displaystyle \partial 
x}=0,\label{MarchesIncomplets4}
\end{equation} \ \par
\ \par

\begin{displaymath} 
\frac{\partial u_{3}}{\partial t}\ -\frac{k^{2}x^{2}}{2}\frac{\partial 
^{2}u_{3}}{\partial x^{2}}-\frac{x^{2}y^{2}}{2}\frac{\partial ^{2}u_{3}}{\partial 
y^{2}}-\rho kx^{2}y\frac{\partial ^{2}u_{3}}{\partial x\partial y}-g(x)\frac{\partial 
u_{3}}{\partial x}-xf(x)y\frac{\partial u_{3}}{\partial y}\end{displaymath} 
\ \par

\begin{equation} 
-(1-\rho ^{2})k^{2}x^{2}exp(u_{1})\displaystyle \left({\displaystyle \frac{\displaystyle 
\partial u_{2}}{\displaystyle \partial x}}\right) ^{2}=0,\label{MarchesIncomplets5}
\end{equation} \ \par
with the initial conditions: \ \par
$u_{1}(0)=0;\ u_{2}(0)=F(x,y);\ u_{3}(0)=0$.( $F$ is the payoff function).\ 
\par
The outline of the paper is as follows:\ \par
{\hskip 1.8em}In section 2, we solve ~(\ref{MarchesIncomplets0}). The 
different derivative terms will be treated separately in order to obtain 
the $L^{\infty }$- stability of the scheme. We prove the convergence 
 of the numerical solution towards a weak solution of the problem. Besides 
the uniqueness of this weak solution is obtained.\ \par
{\hskip 1.8em}In sections 3 and 4, we study ~(\ref{MarchesIncomplets4}), 
~(\ref{MarchesIncomplets5}). We use a change of unknown which lead to 
a variationnel formulation and obtain the existence of a unique solution 
in suitable weighted Sobolev spaces. These  equations are discretized 
by using a backward Euler method in time and a finite element method 
in space. Numerical results are presented.\ \par

\section{Computation of $u_{1}$}
\setcounter{equation}{0}\subsection{Definition of the numerical solution{\hskip 
1.8em}}
In order to solve ~(\ref{MarchesIncomplets0}), we use suitable weighted 
Sobolev spaces, such that no boundary condition is needed in $0$ and 
the function $u_{1}$ has the correct behaviour at infinity.\ \par
To simplify  notation, we denote:\ \par
$\displaystyle F_{1}(x)=f^{2}(x),\ x>0$, $\displaystyle \lambda =k^{2}(\rho ^{2}-\frac{\displaystyle 1}{\displaystyle 2}),$\ 
\par
this coefficient may be positive or negative since $\rho $ is a correlation 
factor and then $\rho $ lies in $[-1,+1]$.\ \par
The equation ~(\ref{MarchesIncomplets0}) becomes:\ \par

\begin{equation} 
\frac{\displaystyle \partial u_{1}}{\displaystyle \partial t}-\frac{\displaystyle 
1}{\displaystyle 2}k^{2}x^{2}\frac{\displaystyle \partial ^{2}u_{1}}{\displaystyle 
\partial x^{2}}-g_{1}(x)\frac{\displaystyle \partial u_{1}}{\displaystyle 
\partial x}+\lambda x^{2}\displaystyle \left({\displaystyle \frac{\displaystyle 
\partial u_{1}}{\displaystyle \partial x}}\right) ^{2}+F_{1}(x)=0.\label{MarchesIncomplets6}
\end{equation} \ \par
\ \par
We will make the following assumptions on the functions $f$  and $g_{1}$:\ 
\par
{\hskip 1.8em}-$f\in W^{1,\infty }({\mathbb{R}}^{+}),\ xf'\in L^{\infty }({\mathbb{R}}^{+}).$\ 
\par
{\hskip 1.8em}-$F_{1}$ is a nonincreasing function, $F_{1}'$ is a bounded 
variation function;  if we denote $\hat F$  the function defined by $\hat F(x)=xF'_{1}(x),\ x>0$,we 
get  $\hat F\in L^{\infty }({\mathbb{R}}^{+})$.\ \par
{\hskip 1.8em}-$g_{1}$ may be written $g_{1}(x)=x\phi (x)$ with $\displaystyle 
\phi (x)=-\delta (x-\sigma _{1})-2\rho kf(x),\ \delta >0$; so, there 
exists $\sigma _{2}>0$ such that $\phi $ is  negative and nonincreasing 
on $[\sigma _{2},+\infty [$ and bounded on $[0,\sigma _{2}]$. \ \par
\ \par
We define the two constants $c_{1}$ and $c_{2}$ by 
\begin{equation} 
c_{1}=\underset{\displaystyle x\in [0,\sigma _{2}]}{\displaystyle {\mathrm{sup}}}\displaystyle 
\left\vert{\displaystyle g'_{1}(x)}\right\vert ,\ c_{2}=\underset{\displaystyle 
x\in [0,\sigma _{2}]}{\displaystyle {\mathrm{sup}}}\displaystyle \left\vert{\displaystyle 
\phi (x)}\right\vert \label{MarchesIncomplets33}
\end{equation}  \ \par
{\hskip 1.8em}We denote by $\Delta t_{n}$ the time increment between 
the levels $t_{n}$ and $t_{n+1}$, $n\geq 0$ and by $u_{1h}^{n}$ the approximate 
solution at the time level $t_{n}$. This solution will be in a finite-dimensional 
space $V_{1h}$ which will be defined below.\ \par
The solution $u_{1h}^{n+1}$ at the time level $t_{n+1}$ is computed in 
two steps: knowing $u_{1h}^{n}$, we compute $u_{1h}^{n+\frac{1}{2}}$, 
approximate solution of 
\begin{equation} 
\frac{\displaystyle \partial u_{1}}{\displaystyle \partial t}+\lambda 
x^{2}\displaystyle \left({\displaystyle \frac{\displaystyle \partial 
u_{1}}{\displaystyle \partial x}}\right) ^{2}=0\label{MarchesIncomplets9}
\end{equation}  obtained by using an explicit upwind scheme. Then starting 
with this intermediate value, we use a backward Euler method in time 
to compute $u_{1h}^{n+1}$; the second order term in ~(\ref{MarchesIncomplets6}) 
is discretized by using a $P_{1}$-finite element method ~\cite{Ciarlet} 
 and the linear first order term by an implicit upwind scheme ~\cite{Ritchmyer} 
in order to get the $L^{\infty }$ -stability of the scheme.\ \par
\ \par
In order to define the finite-dimensional space $V_{1h}$, we first study 
the parabolic problem:\ \par

\begin{equation} 
\frac{\displaystyle \partial u_{1}}{\displaystyle \partial t}-\frac{\displaystyle 
1}{\displaystyle 2}k^{2}x^{2}\frac{\displaystyle \partial ^{2}u_{1}}{\displaystyle 
\partial x^{2}}+F_{1}(x)=0\label{MarchesIncomplets7}
\end{equation} \ \par
to obtain a variational formulation in weighted Sobolev spaces.\ \par
\ \par
\subsubsection{Variational formulation of ~(\ref{MarchesIncomplets7})}
\ \par
Let us consider the two spaces:\ \par

\begin{displaymath} 
H_{1}=\left\lbrace{\left.{v\in {\mathcal{D}}'({\mathbb{R}}^{+})/\frac{v}{1+x}\in 
L^{2}({\mathbb{R}}^{+})}\right\rbrace }\right. \end{displaymath} \ \par
and  
\begin{displaymath} 
V_{1}=\left\lbrace{\left.{v\in {\mathcal{D}}'({\mathbb{R}}^{+})/\frac{v}{1+x}\in 
L^{2}({\mathbb{R}}^{+}),\ \frac{xv'}{1+x}\in L^{2}({\mathbb{R}}^{+})}\right\rbrace 
}\right. .\end{displaymath} \ \par
The space $H_{1}$ is equipped with the following scalar product:\ \par
$\displaystyle \forall v,w\in H_{1},\ (v,w)_{1}=\displaystyle \int_{0}^{+\infty }{\frac{\displaystyle 
v(x)w(x)}{\displaystyle (1+x)^{2}}dx}$  and  the associated norm.\ \par
The space $V_{1}$  with  the norm \ \par

\begin{displaymath} 
\left\Vert{v}\right\Vert _{V_{1}}=\left({\int_{0}^{\infty }{\left({\frac{v^{2}(x)}{(1+x)^{2}}+\frac{x^{2}}{(1+x)^{2}}\left({\frac{dv}{dx}}\right) 
^{2}}\right) dx}}\right) ^{\frac{1}{2}}\end{displaymath} \ \par
is a Hilbert space and$\ {\mathcal{D}}({\mathbb{R}}^{+})$ is dense in 
$V_{1}$  ~\cite{Bolley}{\it{.}}\ \par
We define on $V_{1}\times V_{1}$ the bilinear form: \ \par

\begin{displaymath} 
\forall v,w\in V_{1},\ a(v,w)=\frac{1}{2}k^{2}\int_{0}^{+\infty }{\frac{dv}{dx}\frac{d}{dx}\left({w\frac{x^{2}}{(1+x)^{2}}}\right) 
dx}\end{displaymath} \ \par
or 
\begin{displaymath} 
a(v,w)=\frac{1}{2}k^{2}\int_{0}^{+\infty }{\frac{x^{2}}{(1+x)^{2}}\frac{dv}{dx}\frac{dw}{dx}dx}+k^{2}\int_{0}^{+\infty 
}{\frac{x}{(1+x)^{3}}\frac{dv}{dx}wdx}.\end{displaymath} \ \par
This bilinear form is continue on $V_{1}\times V_{1}$ and we have the 
equality\ \par
$\displaystyle a(v,v)=\frac{\displaystyle k^{2}}{\displaystyle 2}\displaystyle \int_{0}^{+\infty 
}{\frac{\displaystyle x^{2}}{\displaystyle (1+x)^{2}}\displaystyle \left({\displaystyle 
\frac{\displaystyle dv}{\displaystyle dx}}\right) ^{2}dx}-\frac{\displaystyle 
k^{2}}{\displaystyle 2}\displaystyle \int_{0}^{+\infty }{\frac{\displaystyle 
1-2x}{\displaystyle (1+x)^{4}}v^{2}dx},$   \ \par
then we get: 
\begin{displaymath} 
\forall v\in V_{1},\ a(v,v)\geq \frac{k^{2}}{2}\left\Vert{v}\right\Vert 
_{V_{1}}^{2}-k^{2}\left\Vert{v}\right\Vert _{H_{1}}^{2}.\end{displaymath} 
\ \par
Since the function $F_{1}$ is in $H_{1}$, the following variational problem:\ 
\par
Find $\displaystyle u_{1}\in L^{2}(0,T;V_{1})\cap C(0,T;H_{1})$ such 
that:\ \par

\begin{equation} 
\displaystyle \left\lbrace{\displaystyle \begin{matrix}{\displaystyle 
\displaystyle \left({\displaystyle \frac{\displaystyle \partial u_{1}}{\displaystyle 
\partial t},v}\right) _{1}+a(u_{1},v)=-(F_{1},v)_{1},\ \forall v\in V_{1}}\cr 
{\displaystyle u_{1}(0)=0}\cr \end{matrix}.}\right. \label{MarchesIncomplets8}
\end{equation} \ \par
has a unique solution ~\cite{Lions}.\ \par
\ \par
\subsubsection{Approximation of ~(\ref{MarchesIncomplets8})}
\ \par
{\hskip 1.8em}The finite-dimensional space $V_{1h}$ will be a subspace 
of $V_{1}$ defined in the following way:\ \par
Let $(x_{i})_{0\leq i\leq N}$ an increasing sequence ($x_{0}=0)$. We 
denote $h_{i}=x_{i}-x_{i-1}$, $I_{i}=(x_{i-1},x_{i}),\ 1\leq i\leq N$, 
$I_{N+1}=(x_{N},+\infty )$\ \par

\begin{displaymath} 
V_{1h}=\left\lbrace{\left.{v_{h}\in C^{0}({\mathbb{R}}^{+})/\ \ v_{h\mid 
I_{i}}\in P_{1},\ 1\leq i\leq N,\ v_{h\mid I_{N+1}}\in P_{0}}\right\rbrace 
}\right. .\end{displaymath} \ \par
The variable $x$  is the volatility  which lies, in practice, in $]0,1[$, 
so, we may use a constant space step $h$ on $(0,1)$ and  an increasing 
sequence ($h_{i})\ $ for $x\geq 1$ in order that the number of nodes 
is not too important.\ \par
If $v_{h}\in V_{1h}$, we denote $v_{i}=v_{h}(x_{i})$.\ \par
We define on $V_{1h}\ $an approximate scalar product:\ \par

\begin{displaymath} 
(v_{h},w_{h})_{h}=\ \frac{h_{1}}{2}v_{0}w_{0}+\sum_{i=1}^{N-1}{\frac{h_{i}+h_{i+1}}{2}\frac{1}{(1+x_{i})^{2}}v_{i}w_{i}}+v_{N}w_{N}\left({\frac{h_{N}}{2(1+x_{N})^{2}}+\frac{1}{1+x_{N}}}\right) 
\end{displaymath} \ \par
obtained by using the trapezoid method on each interval $I_{i}$, $1\leq i\leq N$ 
; the last integral being computed exactly.\ \par
We also define the Lagrange interpolate $\pi _{h}F_{1}$ of $F_{1}$ by:\ 
\par
$\pi _{h}F_{1}\in V_{1h}$ and $\pi _{h}F_{1}(x_{i})=F_{1}(x_{i})=F_{1i},$ 
$\ 1\leq i\leq N$.\ \par
The approximate solution of ~(\ref{MarchesIncomplets8}) at the time level 
$t_{n+1}$ is the solution of :\ \par

\begin{displaymath} 
(u_{1h}^{n+1},v_{h})_{h}+\Delta t_{n}a(u_{1h}^{n+1},v_{h})=(u_{1h}^{n},v_{h})_{h}-\Delta 
t_{n}(\pi _{h}F_{1},v_{h})_{h},\ \forall v_{h}\in V_{1h},\end{displaymath} 
\ \par

\begin{displaymath} 
u_{1h}^{0}=0.\end{displaymath} \ \par
This may be written:
\begin{displaymath} 
u_{10}^{n+1}=u_{10}^{n}-\Delta t_{n}F_{10},\end{displaymath} \ \par

\begin{displaymath} 
u_{1i}^{n+1}+\Delta t_{n}\alpha _{i}\left({\left({\frac{1}{h_{i}}+\frac{1}{h_{i+1}}}\right) 
u_{1i}^{n+1}-\frac{1}{h_{i}}u_{1,i-1}^{n+1}-\frac{1}{h_{i+1}}u_{1,i+1}^{n+1}}\right) 
=u_{1i}^{n}-\Delta t_{n}F_{1i},\end{displaymath} \ \par
\ \par

\begin{displaymath} 
u_{1N}^{n+1}+\frac{\Delta t_{n}}{h_{N}}\alpha _{N}\left({u_{1N}^{n+1}-u_{1,N-1}^{n+1}}\right) 
=u_{1N}^{n}-\Delta t_{n}F_{1N},\end{displaymath} \ \par
with $\displaystyle \alpha _{i}=\frac{\displaystyle k^{2}x_{i}^{2}}{\displaystyle h_{i}+h_{i+1}},\ 
1\leq i\leq N-1,\ \alpha _{N}=\frac{\displaystyle k^{2}x_{N}^{2}}{\displaystyle 
h_{N}+2(1+x_{N})}$.\ \par
\ \par
\subsubsection{Approximation of the first order terms}
We compute now an approximate solution of ~(\ref{MarchesIncomplets9}) 
by using an explicit upwind scheme.\ \par
Let us denote by $v_{h}^{n}$ the derivative of $u_{1h}^{n}$\ \par
$\displaystyle v^{n}_{\displaystyle h\mid I_{i}}=v_{i}^{n}=\frac{\displaystyle u_{1i}^{n}-u_{1,i-1}^{n}}{\displaystyle 
h_{i}}\ ,\ \ 1\leq i\leq N,\ v^{n}_{\displaystyle h\mid I_{N+1}}=v_{N+1}^{n}=0$.\ 
\par
We set $v_{0}^{n}=0.$\ \par
We shall prove below that the function $v_{h}^{n}$ is positive and the 
function $u_{1h}^{n}$ is negative; we define $u_{1h}^{n+\frac{1}{2}}\in V_{1h}$ 
by:
\begin{equation} 
u_{1i}^{n+\frac{1}{2}}=u_{1i}^{n}-\lambda \Delta t_{n}x_{i}^{2}(v_{i}^{n})^{2}\ 
\ if\ \ \lambda >0,\label{MarchesIncomplets11}
\end{equation} \ \par

\begin{equation} 
u_{1i}^{n+\frac{1}{2}}=u_{1i}^{n}-\lambda \Delta t_{n}x_{i+1}^{2}(v_{i+1}^{n})^{2}\ 
if\ \lambda <0,\label{MarchesIncomplets12}
\end{equation} \ \par
$0\leq i\leq N$. \ \par
For the linear first order term, since the function $g_{1}$ is not bounded, 
we use an implicit scheme, which will be decentered in order to get a 
monotone matrix.\ \par
\ \par
Finally, the solution $u_{1h}^{n+1}\in V_{1h}$ of ~(\ref{MarchesIncomplets6}) 
 is defined by:
\begin{displaymath} 
u_{10}^{n+1}=u_{10}^{n+\frac{1}{2}}-\Delta t_{n}F_{10},\end{displaymath} 
\ \par
\ \par

\begin{displaymath} 
u_{1i}^{n+1}+\Delta t_{n}\left({\alpha _{i}\left({\frac{1}{h_{i}}+\frac{1}{h_{i+1}}}\right) 
+\gamma _{i}\left({\frac{1-\delta _{i}}{h_{i+1}}-\frac{\delta _{i}}{h_{i}}}\right) 
}\right) u_{1i}^{n+1}\end{displaymath}  $\displaystyle \ \ \ \ \ \ \ \ \ \ \ \ \ \ \ \ \ \ \ \ \ \ \ \ \ \ \ \ \ \ -\frac{\displaystyle 
\Delta t_{n}}{\displaystyle h_{i}}\displaystyle \left({\displaystyle 
\alpha _{i}-\gamma _{i}\delta i}\right) u_{1,i-1}^{n+1}-\frac{\displaystyle 
\Delta t_{n}}{\displaystyle h_{i+1}}\displaystyle \left({\displaystyle 
\alpha _{i}+\gamma _{i}(1-\delta _{i})}\right) u_{1,i+1}^{n+1}$      
                      \ \par

\begin{equation} 
=u_{1i}^{n+\frac{1}{2}}-\Delta t_{n}F_{1i},\label{MarchesIncomplets10}
\end{equation} \ \par
\ \par

\begin{displaymath} 
u_{1N}^{n+1}+\frac{\Delta t_{n}}{h_{N}}\left({\alpha _{N}-\gamma _{N}}\right) 
(u_{1N}^{n+1}-u_{1,N-1}^{n+1})\end{displaymath}  \ \par

\begin{displaymath} 
=u_{1N}^{n+\frac{1}{2}}-\Delta t_{n}F_{1N},\end{displaymath} \ \par
with $\gamma _{i}=g_{1}(x_{i}),\ 1\leq i\leq N,$\ \par
\ \par
$\displaystyle \delta _{i}=\displaystyle \left\lbrace{\displaystyle \begin{matrix}{\displaystyle 
0\ if\ \gamma _{i}\geq 0,}\cr {\displaystyle 1\ if\ \gamma _{i}<0.}\cr 
\end{matrix}}\right. $\ \par
\ \par
\ \par
Since $g_{1}$ is negative for $x\geq \sigma _{2}$ , we get  $\delta _{i}=1$ 
if $i$ is large enough.   \ \par
The preceding equations may be written by using the derivative $v_{h}^{n}$:\ 
\par

\begin{displaymath} 
u_{1i}^{n+1}+\Delta t_{n}\left({\left({\alpha _{i}-\delta _{i}\gamma _{i}}\right) 
v_{i}^{n+1}-\left({\alpha _{i}+(1-\delta _{i})\gamma _{i}}\right) v_{i+1}^{n+1}}\right) 
=u_{1i}^{n+\frac{1}{2}}-\Delta t_{n}F_{1i},\end{displaymath} \ \par
$0\leq i\leq N$\ \par
\ \par
\subsection{Properties of the scheme}
We prove that under a stability condition, the approximate solution $u_{1h}^{n}$ 
is negative and its derivative $v_{h}^{n}$ is positive.\ \par
Let us denote $U_{1h}^{n}\ $ the vector of ${\mathbb{R}}^{N+1}$ of components 
$(u_{1i}^{n}),\ 0\leq i\leq N\ $and $F_{1h}$ the vector of components 
$(F_{1i}),\ 0\leq i\leq N$.\ \par
The numerical scheme ~(\ref{MarchesIncomplets10}) may be written:
\begin{equation} 
(I+\Delta t_{n}A_{h})U_{1h}^{n+1}=U_{1h}^{n+\frac{1}{2}}-\Delta t_{n}F_{1h}\label{MarchesIncomplets19}
\end{equation} \ \par
where the matrix $A_{h}$ is tridiagonal and monotone.\ \par
\ \par
From  ~(\ref{MarchesIncomplets11}), ~(\ref{MarchesIncomplets12}) and 
~(\ref{MarchesIncomplets10}), we get immediately the equations satisfied 
by $v_{h}^{n}$:\ \par
If $v_{h}^{n+\frac{1}{2}}$ denotes the derivative of $u_{h}^{n+\frac{1}{2}}$, 
we obtain for $1\leq i\leq N$\ \par

\begin{equation} 
v_{i}^{n+\frac{1}{2}}=v_{i}^{n}-\lambda \frac{\Delta t_{n}}{h_{i}}\left({x_{i}^{2}(v_{i}^{n})^{2}-x_{i-1}^{2}(v_{i-1}^{n})^{2}}\right) 
\ \ if\ \lambda >0,\label{MarchesIncomplets13}
\end{equation} \ \par

\begin{equation} 
v_{i}^{n+\frac{1}{2}}=v_{i}^{n}-\lambda \frac{\Delta t_{n}}{h_{i}}\left({x_{i+1}^{2}(v_{i+1}^{n})^{2}-x_{i}^{2}(v_{i}^{n})^{2}}\right) 
\ \ if\ \lambda <0,\label{MarchesIncomplets14}
\end{equation} \ \par
and $v_{h}^{n+1}$ satisfies: 
\begin{displaymath} 
v_{i}^{n+1}+\frac{\Delta t_{n}}{h_{i}}\left({\alpha _{i}+\alpha _{i-1}-\delta 
_{i}\gamma _{i}+(1-\delta _{i-1})\gamma _{i-1}}\right) v_{i}^{n+1}\end{displaymath} 
\ \par

\begin{displaymath} 
-\frac{\displaystyle \Delta t_{n}}{\displaystyle h_{i}}\displaystyle \left({\displaystyle 
\alpha _{i-1}-\delta _{i-1}\gamma _{i-1}}\right) v_{i-1}^{n+1}-\frac{\displaystyle 
\Delta t_{n}}{\displaystyle h_{i}}\displaystyle \left({\displaystyle 
\alpha _{i}+(1-\delta _{i})\gamma _{i}}\right) v_{i+1}^{n+1}\end{displaymath} 
\ \par

\begin{equation} 
=v_{i}^{n+\frac{1}{2}}-\frac{\Delta t_{n}}{h_{i}}\left({F_{1i}-F_{1},_{i-1}}\right) 
\label{MarchesIncomplets21}
\end{equation} \ \par
for $1\leq i\leq N$\ \par
which may be written: 
\begin{equation} 
(I+\Delta t_{n}B_{h})V_{h}^{n+1}=V_{h}^{n+\frac{1}{2}}-\Delta t_{n}G_{1h},\label{MarchesIncomplets37}
\end{equation} \ \par
\ \par
where $V_{h}^{n}$ is the vector of ${\mathbb{R}}^{N}$ of components $(v_{i}^{n}),\ 1\leq i\leq N,$ 
$B_{h}$ is a tridiagonal matrix $(N\times N),$\ \par
$G_{1h}$ is the vector of ${\mathbb{R}}^{N}$ of components$\displaystyle 
\ \frac{\displaystyle F_{1i}-F_{1,i-1}}{\displaystyle h_{i}},\ 1\leq i\leq 
N.$\ \par
\ \par
\begin{Proposition} \label{MarchesIncomplets17}If the following stability 
condition\ \par

\begin{equation} 
\underset{\displaystyle i\geq 1}{\displaystyle {\mathrm{sup}}}\displaystyle 
\left\vert{\displaystyle \lambda \frac{\displaystyle \Delta t_{n}}{\displaystyle 
h_{i}}x_{i}^{2}v_{i}^{n}}\right\vert \leq 1\ \label{MarchesIncomplets15}
\end{equation} \ \par
is satisfied and  
\begin{equation} 
\Delta t_{n}c_{1}<1,\label{MarchesIncomplets16}
\end{equation}  then the function $v_{h}^{n}$ is nonnegative for $n\geq 0$. 
\ \par
\end{Proposition}
\ \par
{\it{Proof: }} We can rewrite ~(\ref{MarchesIncomplets13}) as:\ \par

\begin{equation} 
v_{i}^{n+\frac{1}{2}}=v_{i}^{n}\left({1-\lambda \frac{\Delta t_{n}}{h_{i}}x_{i}^{2}v_{i}^{n}}\right) 
+\lambda \frac{\Delta t_{n}}{h_{i}}x_{i-1}^{2}(v_{i-1}^{n})^{2}\ si\ 
\lambda >0\label{MarchesIncomplets20}
\end{equation} \ \par
and we have an analogous formula  for  $\lambda <0.$\ \par
If ~(\ref{MarchesIncomplets15}) is satisfied, we get immediately that 
$v_{h}^{n}\geq 0$ implies $v_{h}^{n+\frac{1}{2}}\geq 0.$Since $F_{1}$ 
is  decreasing, the vector $V_{h}^{n+\frac{1}{2}}-\Delta t_{n}G_{1h}$ 
is nonnegative and the vector $V_{h}^{n+1}$ will be nonnegative if  $I+\Delta t_{n}B_{h}$ 
is a monotone matrix; this will be true if $1-\frac{\Delta t_{n}}{h_{i}}(\gamma _{i}-\gamma _{i-1})>0$, 
for $1\leq i\leq N$. Since $g_{1}$ is a decreasing function for $x\geq \sigma _{2},\ $this 
condition is satisfied for $i$ large enough and if $x\leq \sigma _{2}$,, 
we get $1-\frac{\Delta t_{n}}{h_{i}}(\gamma _{i}-\gamma _{i-1})\geq 1-c_{1}\Delta 
t_{n}>0$ from ~(\ref{MarchesIncomplets16}). \ \par
\ \par
We deduce immediately the following results:\ \par
\ \par
\begin{Proposition}  Under the hypotheses of  proposition ~\ref{MarchesIncomplets17}, 
the numerical solution $u_{1h}^{n}$ satisfies: $u_{1h}^{n}\leq 0$ for 
$n\geq 0$. \ \par
\end{Proposition}
\ \par
{\it{Proof: }}We can rewrite ~(\ref{MarchesIncomplets11}) as \ \par

\begin{equation} 
u_{1i}^{n+\frac{1}{2}}=u_{1i}^{n}\left({1-\lambda \frac{\Delta t_{n}}{h_{i}}x_{i}^{2}v_{i}^{n}}\right) 
+\lambda \frac{\Delta t_{n}}{h_{i}}x_{i}^{2}v_{i}^{n}u_{1,i-1}^{n}\ if\ 
\lambda >0\label{MarchesIncomplets18}
\end{equation} \ \par
and we have an analogous equality for $\lambda <0$. From proposition 
~(\ref{MarchesIncomplets17}), the function $v_{h}^{n}$ is nonnegative, 
hence if $u_{1h}^{n}\leq 0$, we get $u_{1h}^{n+\frac{1}{2}}\leq 0$.\ 
\par
Since $I+\Delta t_{n}A_{h}$ is a monotone matrix and $F_{1}\geq 0$, we 
deduce that $u_{1h}^{n+1}\leq 0$. \ \par
\ \par
\begin{Proposition} Under the hypotheses of proposition ~\ref{MarchesIncomplets17}, 
the numerical solution satisfies\ \par

\begin{displaymath} 
\displaystyle \left\Vert{\displaystyle u_{1h}^{n}}\right\Vert _{\displaystyle 
L^{\infty }({\mathbb{R}}^{+})}\leq t_{n}F_{1}(0)\end{displaymath}  for 
$n\geq 0.$\ \par
\end{Proposition}
\ \par
{\it{Proof: }} We get immediately from ~(\ref{MarchesIncomplets18}): 
$\left\Vert{u_{1h}^{n+\frac{1}{2}}}\right\Vert _{L^{\infty }({\mathbb{R}}^{+})}\leq 
\left\Vert{u_{1h}^{n}}\right\Vert _{L^{\infty }({\mathbb{R}}^{+})}$ and 
since $A_{h}$ is a monotone matrix, it follows from ~(\ref{MarchesIncomplets19}) 
that\ \par
$\displaystyle \displaystyle \left\Vert{\displaystyle u_{1h}^{n+1}}\right\Vert _{\displaystyle 
L^{\infty }({\mathbb{R}}^{+})}$ $\leq \left\Vert{u_{1h}^{n+\frac{1}{2}}}\right\Vert _{L^{\infty }({\mathbb{R}}^{+})}+\Delta 
t_{n}F_{1}(0)$  which completes the proof.\ \par
\ \par
\ \par
\begin{Proposition} Under the hypotheses of proposition  ~\ref{MarchesIncomplets17}, 
the function $v_{h}^{n}$ satisfies:\ \par

\begin{displaymath} 
\displaystyle \left\Vert{\displaystyle v_{h}^{n}}\right\Vert _{\displaystyle 
L^{1}({\mathbb{R}}^{+})}\leq t_{n}F_{1}(0)\end{displaymath}  for $n\geq 0$.\ 
\par
\end{Proposition}
\ \par
{\it{Proof: }}From ~(\ref{MarchesIncomplets20}), we get if $\lambda >0$ 
\ \par

\begin{displaymath} 
\sum_{i=1}^{N}{h_{i}v_{i}^{n+\frac{1}{2}}}=\sum_{i=1}^{N}{h_{i}v_{i}^{n}\left({1-\lambda 
\frac{\Delta t_{n}}{h_{i}}x_{i}^{2}v_{i}^{n}}\right) }+\lambda \Delta 
t_{n}\sum_{i=1}^{N}{x_{i-1}^{2}(v_{i-1}^{n})^{2}},\end{displaymath} \ 
\par
hence $\ \ \ \left\Vert{v_{h}^{n+\frac{1}{2}}}\right\Vert _{L^{1}({\mathbb{R}}^{+})}\leq 
\left\Vert{v_{h}^{n}}\right\Vert _{L^{1}({\mathbb{R}}^{+})}.$\ \par
For $\lambda <0$, we obtain the same inequality.\ \par
\ \par
Besides from ~(\ref{MarchesIncomplets21}), we get: \ \par
$\left\Vert{v_{h}^{n+1}}\right\Vert _{L^{1}({\mathbb{R}}^{+})}\leq \left\Vert{v_{h}^{n+\frac{1}{2}}}\right\Vert 
_{L^{1}({\mathbb{R}}^{+})}-\Delta t_{n}F_{1N}+\Delta t_{n}F_{10}$.\ \par
and we deduce the result.\ \par
\ \par
We prove now that under some hypothesis on the sequence $(h_{i})$, the 
function $xv_{h}^{n}$ is bounded in $L^{\infty }({\mathbb{R}}^{+})$ and 
it is possible to choose the nodes $(x_{i})_{1\leq i\leq N}$ such that 
the stability condition is not too restrictive.\ \par
\ \par
We define the function $\hat v_{h}^{n}$ by:\ \par
$\displaystyle \hat v^{n}_{\displaystyle h\mid I_{i}}=\hat v_{i}^{n}=x_{i}v_{i}^{n},\ 
1\leq i\leq N,\hat v^{n}_{\displaystyle h\mid I_{N+1}}=\hat v_{N+1}^{n}=0$.\ 
\par
\ \par
\begin{Proposition} \label{MarchesIncomplets36}If the following stability 
condition 
\begin{equation} 
\underset{\displaystyle i\geq 1}{\displaystyle {\mathrm{sup}}}\lambda 
\frac{\displaystyle \Delta t_{n}}{\displaystyle h_{i}}x_{i}(\hat v_{i}^{n}+\hat 
v_{i-1}^{n})\leq 1\ if\ \ \lambda >0\ and\ \underset{\displaystyle i\geq 
1}{\displaystyle {\mathrm{sup}}}\displaystyle \left\vert{\displaystyle 
\lambda }\right\vert \frac{\displaystyle \Delta t_{n}}{\displaystyle 
h_{i}}x_{i}(\hat v_{i}^{n}+\hat v_{i+1}^{n})\leq 1\ if\ \ \lambda <0\label{MarchesIncomplets28}
\end{equation}  is satisfied and if the sequence $(h_{i})_{1\leq i\leq N}$ 
satisfy: There exists a positive constant $c$ such that 
\begin{equation} 
\frac{\displaystyle x_{i}}{\displaystyle x_{i+1}}\frac{\displaystyle h_{i+1}}{\displaystyle 
h_{i}+h_{i+1}}-\frac{\displaystyle x_{i-1}}{\displaystyle x_{i}}\frac{\displaystyle 
h_{i}}{\displaystyle h_{i}+h_{i-1}}\geq -c\frac{\displaystyle h_{i}}{\displaystyle 
x_{i}},\ 1\leq i\leq N-1\label{MarchesIncomplets30}
\end{equation}  and 
\begin{equation} 
\frac{\displaystyle x_{N}}{\displaystyle h_{N}+2(1+x_{N})}-\frac{\displaystyle 
x_{N-1}}{\displaystyle x_{N}}\frac{\displaystyle h_{N}}{\displaystyle 
h_{N}+h_{N-1}}\geq -c\frac{\displaystyle h_{N}}{\displaystyle x_{N}},\label{MarchesIncomplets31}
\end{equation}  then if $\Delta t\leq \Delta t_{0}$, $\Delta t_{0}$ depending 
on $c,c_{1},c_{2},$ the following estimate holds: 
\begin{equation} 
\displaystyle \left\Vert{\displaystyle \hat v_{h}^{n}}\right\Vert _{\displaystyle 
L^{\infty }({\mathbb{R}}^{+})}\leq e^{\displaystyle Ct_{n}}\displaystyle 
\left\Vert{\displaystyle \hat F}\right\Vert _{\displaystyle L^{\infty 
}({\mathbb{R}}^{+})}\label{MarchesIncomplets34}
\end{equation}  for $n\geq 0$ and $C\ $is a constant depending on $c,c_{1},c_{2}$. 
\ \par
\end{Proposition}
\ \par
{\it{Proof: }}We get immediately from ~(\ref{MarchesIncomplets13})\ \par
$\hat v_{i}^{n+\frac{1}{2}}$ $\displaystyle =\hat v_{i}^{n}-\lambda \frac{\displaystyle \Delta t_{n}}{\displaystyle 
h_{i}}x_{i}\displaystyle \left({\displaystyle (\hat v_{i}^{n})^{2}-(\hat 
v_{i-1}^{n})^{2}}\right) \ if\ \lambda >0$\ \par
which may be written:\ \par

\begin{displaymath} 
\hat v_{i}^{n+\frac{1}{2}}=\hat v_{i}^{n}\left({1-\lambda \frac{\Delta 
t_{n}}{h_{i}}x_{i}(\hat v_{i}^{n}+\hat v_{i-1}^{n})}\right) +\lambda 
\frac{\Delta t_{n}}{h_{i}}x_{i}(\hat v_{i}^{n}+\hat v_{i-1}^{n})\hat 
v_{i-1}^{n}\end{displaymath} \ \par
and if ~(\ref{MarchesIncomplets28}) is satisfied, we obtain:\ \par

\begin{equation} 
\left\Vert{\hat v_{h}^{n+\frac{1}{2}}}\right\Vert _{L^{\infty }({\mathbb{R}}^{+})}\leq 
\left\Vert{\hat v_{h}^{n}}\right\Vert _{L^{\infty }({\mathbb{R}}^{+})}.\label{MarchesIncomplets29}
\end{equation} \ \par
If $\lambda <0,$ we get from ~(\ref{MarchesIncomplets14}) \ \par

\begin{displaymath} 
\hat v_{i}^{n+\frac{1}{2}}=\hat v_{i}^{n}\left({1+\lambda \frac{\Delta 
t_{n}}{h_{i}}x_{i}(\hat v_{i+1}^{n}+\hat v_{i}^{n})}\right) -\lambda 
\frac{\Delta t_{n}}{h_{i}}x_{i}(\hat v_{i+1}^{n}+\hat v_{i}^{n})\hat 
v_{i+1}^{n}\end{displaymath} \ \par
and the estimate ~(\ref{MarchesIncomplets29}) holds if ~(\ref{MarchesIncomplets28}) 
is satisfied.\ \par
\ \par
Further the following equality results of ~(\ref{MarchesIncomplets21}) 
for $1\leq i\leq N$ \ \par

\begin{displaymath} 
\hat v_{i}^{n+1}+\frac{\Delta t_{n}}{h_{i}}\left({\alpha _{i}+\alpha _{i-1}-\delta 
_{i}\gamma _{i}+(1-\delta _{i-1})\gamma _{i-1}}\right) \hat v_{i}^{n+1}\end{displaymath} 
\ \par

\begin{equation} 
-\frac{\displaystyle \Delta t_{n}}{\displaystyle h_{i}}\displaystyle \left({\displaystyle 
\alpha _{i-1}-\delta _{i-1}\gamma _{i-1}}\right) \frac{\displaystyle 
x_{i}}{\displaystyle x_{i-1}}\hat v_{i-1}^{n+1}-\frac{\displaystyle \Delta 
t_{n}}{\displaystyle h_{i}}\displaystyle \left({\displaystyle \alpha 
_{i}+(1-\delta _{i})\gamma _{i}}\right) \frac{\displaystyle x_{i}}{\displaystyle 
x_{i+1}}\hat v_{i+1}^{n+1}
\end{equation}  \ \par

\begin{displaymath} 
=\hat v_{i}^{n+\frac{1}{2}}-\frac{\Delta t_{n}}{h_{i}}x_{i}\left({F_{1i}-F_{1},_{i-1}}\right) 
\label{MarchesIncomplets32}\end{displaymath} \ \par
\ \par
which may be written
\begin{displaymath} 
(I+\Delta t_{n}C_{h})\hat V_{h}^{n+1}=\hat V_{h}^{n+\frac{1}{2}}-\Delta 
t_{n}\hat F_{h}\end{displaymath} \ \par
where  $\hat F_{h}\ $is the  vector of components : $\displaystyle \hat F_{i}=\ x_{i}\ \frac{\displaystyle F_{1i}-F_{1},_{i-1}}{\displaystyle 
h_{i}}$ , $1\leq i\leq N$. $C_{h}$ is a tridiagonal matrix $(N\times N).$\ 
\par
Besides, we have  $c_{ii}>0,\ c_{ij}\leq 0\ for\ i{\not =}j$, $1\leq i\leq N$.\ 
\par
\ \par
We prove that if ~(\ref{MarchesIncomplets30}) and ~(\ref{MarchesIncomplets31}) 
are satisfied, there exists a positive constant $\hat c$ depending on 
$c_{1},\ c_{2},\ c$ such that: $\displaystyle \displaystyle \sum_{j=1}^{N}{c_{ij}\geq -\hat c},\ 1\leq i\leq N$.\ 
\par
We deduce from ~(\ref{MarchesIncomplets32}): \ \par
$\displaystyle \displaystyle \sum_{j=1}^{N}{c_{ij}}=\frac{\displaystyle 1}{\displaystyle 
h_{i}}\displaystyle \left({\displaystyle \alpha _{i}\frac{\displaystyle 
h_{i+1}}{\displaystyle x_{i+1}}-\alpha _{i-1}\frac{\displaystyle h_{i}}{\displaystyle 
x_{i-1}}-\gamma _{i}+\gamma _{i-1}+\gamma _{i}(1-\delta _{i})\frac{\displaystyle 
h_{i+1}}{\displaystyle x_{i+1}}+\gamma _{i-1}\delta _{i-1}\frac{\displaystyle 
h_{i}}{\displaystyle x_{i-1}}}\right) ,\ 1\leq i\leq N-1$ \ \par
\ \par
$\displaystyle \displaystyle \sum_{j=1}^{N}{c_{Nj}}=\frac{\displaystyle 1}{\displaystyle 
h_{N}}\displaystyle \left({\displaystyle \alpha _{N}-\alpha _{N-1}\frac{\displaystyle 
h_{N}}{\displaystyle x_{N-1}}-\gamma _{N}+\gamma _{N-1}\frac{\displaystyle 
x_{N}}{\displaystyle x_{N-1}}}\right) .$\ \par
\ \par
\ \par
Let us denote  $\displaystyle A_{i}^{1}=\frac{\displaystyle 1}{\displaystyle h_{i}}\displaystyle \left({\displaystyle 
\alpha _{i}\frac{\displaystyle h_{i+1}}{\displaystyle x_{i+1}}-\alpha 
_{i-1}\frac{\displaystyle h_{i}}{\displaystyle x_{i-1}}}\right) =\frac{\displaystyle 
k^{2}x_{i}}{\displaystyle h_{i}}\displaystyle \left({\displaystyle \frac{\displaystyle 
x_{i}}{\displaystyle x_{i+1}}\frac{\displaystyle h_{i+1}}{\displaystyle 
h_{i}+h_{i+1}}-\frac{\displaystyle x_{i-1}}{\displaystyle x_{i}}\frac{\displaystyle 
h_{i}}{\displaystyle h_{i-1}+h_{i}}}\right) $,\ \par
$1\leq i\leq N-1$,\ \par
$\displaystyle A_{N}^{1}=\frac{\displaystyle 1}{\displaystyle h_{N}}\displaystyle \left({\displaystyle 
\alpha _{N}-\alpha _{N-1}\frac{\displaystyle h_{N}}{\displaystyle x_{N-1}}}\right) 
=\frac{\displaystyle k^{2}x_{N}}{\displaystyle h_{N}}\displaystyle \left({\displaystyle 
\frac{\displaystyle x_{N}}{\displaystyle h_{N}+2(1+x_{N})}-\frac{\displaystyle 
x_{N-1}}{\displaystyle x_{N}}\frac{\displaystyle h_{N}}{\displaystyle 
h_{N}+h_{N-1}}}\right) $; \ \par
\ \par
and  $\displaystyle A_{i}^{2}=\frac{\displaystyle 1}{\displaystyle h_{i}}\displaystyle \left({\displaystyle 
-\gamma _{i}+\gamma _{i-1}+\gamma _{i}(1-\delta _{i})\frac{\displaystyle 
h_{i+1}}{\displaystyle x_{i+1}}+\gamma _{i-1}\delta _{i-1}\frac{\displaystyle 
h_{i}}{\displaystyle x_{i-1}}}\right) \ 2\leq i\leq N$, \ \par
$\displaystyle A_{1}^{2}=\frac{\displaystyle 1}{\displaystyle h}\displaystyle \left({\displaystyle 
-\delta _{1}-\frac{\displaystyle 1-\delta _{1}}{\displaystyle 2}}\right) 
\gamma _{1}.\ $\ \par
\ \par
\ \par
We get $\displaystyle \displaystyle \sum_{j=1}^{N}{c_{ij}}=A_{i}^{1}+A_{i}^{2}$.\ 
\par
From ~(\ref{MarchesIncomplets30}) and ~(\ref{MarchesIncomplets31}), it 
follows that $\displaystyle A_{i}^{1}\geq -ck^{2},\ 1\leq i\leq N$.\ 
\par
\ \par
Further , we have :  $\displaystyle A_{1}^{2}=-\phi (x_{1})(\delta _{1}+\frac{\displaystyle 1-\delta _{1}}{\displaystyle 
2})$\ \par
and for $i\geq 2,$  $\displaystyle A_{i}^{2}=-\frac{\displaystyle g_{1}(x_{i})-g_{1}(x_{i-1})}{\displaystyle 
x_{i}-x_{i-1}}+(1-\delta _{i})\frac{\displaystyle x_{i}h_{i+1}}{\displaystyle 
x_{i+1}h_{i}}\phi (x_{i})+\delta _{i-1}\phi (x_{i-1}),\ $\ \par
\ \par
We deduce from ~(\ref{MarchesIncomplets33}) $A_{i}^{2}\geq -(c_{1}+c_{2})\ for\ \ x_{i}\leq \sigma _{2}$.\ 
\par
\ \par
For  $x\geq \sigma _{2}$,  the function $g_{1}$ is negative, so $\delta _{i}=1$, 
 $\displaystyle A_{i}^{2}=\frac{\displaystyle -x_{i}(\phi (x_{i})-\phi (x_{i-1}))}{\displaystyle 
h_{i}}$ \ \par
and  $A_{i}^{2}\geq 0$, since  $\phi $ is decreasing. \ \par
\ \par
Finally, we obtain $\displaystyle \displaystyle \sum_{j=1}^{N}{c_{ij}}\geq -\hat c,\ \ 1\leq i\leq N$, 
with $\displaystyle \hat c=ck^{2}+c_{1}+c_{2}.$\ \par
\ \par
and $\displaystyle \displaystyle \left\Vert{\displaystyle \hat v_{h}^{n+1}}\right\Vert _{\displaystyle 
L^{\infty }({\mathbb{R}}^{+})}\leq \frac{\displaystyle 1}{\displaystyle 
1-\hat c\Delta t_{n}}\ \displaystyle \left({\displaystyle \displaystyle 
\left\Vert{\displaystyle \hat v_{h}^{n}}\right\Vert _{\displaystyle L^{\infty 
}({\mathbb{R}}^{+})}+\Delta t_{n}\displaystyle \left\Vert{\displaystyle 
\hat F}\right\Vert _{\displaystyle L^{\infty }({\mathbb{R}}^{+})}}\right) 
.$\ \par
The estimate ~(\ref{MarchesIncomplets34}) follows.\ \par
\ \par
\ \par
Let us define now  a sequence $(x_{i}),$ satisfying ~(\ref{MarchesIncomplets30}) 
and ~(\ref{MarchesIncomplets31}):\ \par
\ \par
We set: $x_{i}=ih,\ 0\leq i\leq n_{0}$ with $n_{0}h=1$,\ \par
then $x_{i}=e^{\theta h}x_{i-1}$ for $i\geq n_{0}+1$ and $\theta \geq 1$\ 
\par
\ \par
We get: $\displaystyle \frac{\displaystyle x_{i}}{\displaystyle x_{i+1}}\frac{\displaystyle h_{i+1}}{\displaystyle 
h_{i+1}+h_{i}}-\frac{\displaystyle x_{i-1}}{\displaystyle x_{i}}\frac{\displaystyle 
h_{i}}{\displaystyle h_{i}+h_{i-1}}>0,\ i\leq n_{0}+1$, \ \par
\ \par
$\displaystyle \frac{\displaystyle x_{i}}{\displaystyle x_{i+1}}\frac{\displaystyle h_{i+1}}{\displaystyle 
h_{i+1}+h_{i}}-\frac{\displaystyle x_{i-1}}{\displaystyle x_{i}}\frac{\displaystyle 
h_{i}}{\displaystyle h_{i}+h_{i-1}}=0,\ i>n_{0}+1\ $\ \par
\ \par
$\displaystyle \frac{\displaystyle x_{N}}{\displaystyle h_{N}+2(1+x_{N})}-\frac{\displaystyle 
x_{N-1}}{\displaystyle x_{N}}\frac{\displaystyle h_{N}}{\displaystyle 
h_{N}+h_{N-1}}\geq -\frac{\displaystyle 1}{\displaystyle 2x_{N}};\ $ 
then ~(\ref{MarchesIncomplets31}) will be satisfied if $h_{N}\geq \frac{c}{2}$, 
that is $\displaystyle N={\mathcal{O}}\displaystyle \left({\displaystyle \frac{\displaystyle 
\displaystyle \left\vert{\displaystyle {\mathrm{ln}}\ h}\right\vert }{\displaystyle 
h}}\right) $ or $\displaystyle x_{N}={\mathcal{O}}\displaystyle \left({\displaystyle \frac{\displaystyle 
1}{\displaystyle h}}\right) $.\ \par
Besides we have $\displaystyle \frac{\displaystyle x_{i}}{\displaystyle h_{i}}={\mathcal{O}}\displaystyle 
\left({\displaystyle \frac{\displaystyle 1}{\displaystyle h}}\right) 
$  $1\leq i\leq N$, and the stability condition  may be written $\displaystyle 
\frac{\displaystyle \Delta t_{n}}{\displaystyle h}\leq C$, that is the 
classical stability  condition  for hyperbolic problems.\ \par
\ \par
\begin{Proposition} Under the hypotheses of  proposition ~\ref{MarchesIncomplets36}, 
there exists a positive constant $C$ depending on $T,f,g_{1}$ such that 
for $t_{n}\leq T$, the following estimate holds:\ \par

\begin{displaymath} 
\left\Vert{v_{h}^{n}}\right\Vert _{L^{\infty }({\mathbb{R}}^{+})}\leq 
C.\end{displaymath} \ \par
\end{Proposition}
{\it{Proof:}}  For $\lambda >0$, we get from ~(\ref{MarchesIncomplets13}): 
 
\begin{displaymath} 
v_{i}^{n+\frac{1}{2}}=v_{i}^{n}\left({1-\lambda \frac{\Delta t_{n}}{h_{i}}x_{i}(\hat 
v_{i}^{n}+\hat v_{i-1}^{n})}\right) +\lambda \frac{\Delta t_{n}}{h_{i}}x_{i-1}(\hat 
v_{i}^{n}+\hat v_{i-1}^{n})v_{i-1}^{n}\end{displaymath} \ \par
and by using ~(\ref{MarchesIncomplets28}), we obtain:$\left\Vert{v_{h}^{n+\frac{1}{2}}}\right\Vert _{L^{\infty }({\mathbb{R}}^{+})}\leq 
\left\Vert{v_{h}^{n}}\right\Vert _{L^{\infty }({\mathbb{R}}^{+})}.$\ 
\par
For $\lambda <0$, we get:\ \par

\begin{displaymath} 
v_{i}^{n+\frac{1}{2}}=v_{i}^{n}\left({1+\lambda \frac{\Delta t_{n}}{h_{i}}x_{i}(\hat 
v_{i}^{n}+\hat v_{i+1}^{n})}\right) -\lambda \frac{\Delta t_{n}}{h_{i}}x_{i+1}(\hat 
v_{i}^{n}+\hat v_{i+1}^{n})v_{i+1}^{n}\end{displaymath} \ \par
and by using ~(\ref{MarchesIncomplets28}) , we obtain:\ \par
$\left\Vert{v_{h}^{n+\frac{1}{2}}}\right\Vert _{L^{\infty }({\mathbb{R}}^{+})}\leq 
\left\Vert{v_{h}^{n}}\right\Vert _{L^{\infty }({\mathbb{R}}^{+})}\left({1+2\left\vert{\lambda 
}\right\vert \Delta t_{n}\left\Vert{\hat v_{h}^{n}}\right\Vert _{L^{\infty 
}({\mathbb{R}}^{+})}}\right) .$\ \par
\ \par
Besides, it follows from ~(\ref{MarchesIncomplets37}) that\ \par
$(1-c_{1}\Delta t_{n})\left\Vert{v_{h}^{n+1}}\right\Vert _{L^{\infty }({\mathbb{R}}^{+})}\leq 
\left\Vert{v_{h}^{n+\frac{1}{2}}}\right\Vert _{L^{\infty }({\mathbb{R}}^{+})}+\Delta 
t_{n}\left\Vert{F_{1}'}\right\Vert _{L^{\infty }({\mathbb{R}}^{+})}.$\ 
\par
This concludes the proof.\ \par
\ \par
\begin{Proposition} Under the hypotheses of proposition ~\ref{MarchesIncomplets36}, 
there a positive constant $C$ depending on $T,f,g_{1}$ such that for 
$t_{n}\leq T$, the following estimate holds:
\begin{displaymath} 
Var(v_{h}^{n};{\mathbb{R}}^{+})\leq C.\end{displaymath} \ \par
\end{Proposition}
{\it{Proof: }} For $\lambda >0$, we have the equality:\ \par
$v_{i+1}^{n+\frac{1}{2}}-v_{i}^{n+\frac{1}{2}}$ $\displaystyle =v_{i+1}^{n}-v_{i}^{n}-\lambda \frac{\displaystyle \Delta t_{n}}{\displaystyle 
h_{i+1}}\displaystyle \left({\displaystyle x_{i+1}^{2}(v_{i+1}^{n})^{2}-x_{i}^{2}(v_{i}^{n})^{2}}\right) 
+\lambda \frac{\displaystyle \Delta t_{n}}{\displaystyle h_{i}}\displaystyle 
\left({\displaystyle x_{i}^{2}(v_{i}^{n})^{2}-x_{i-1}^{2}(v_{i-1}^{n})^{2}}\right) 
$ \ \par
\ \par
for $1\leq i\leq N-1$, and\ \par
$v_{N+1}^{n+\frac{1}{2}}-v_{N}^{n+\frac{1}{2}}$ $\displaystyle =v_{N+1}^{n}-v_{N}^{n}+\lambda \frac{\displaystyle \Delta t_{n}}{\displaystyle 
h_{N}}\displaystyle \left({\displaystyle x_{N}^{2}(v_{N}^{n})^{2}-x_{N-1}(v_{N-1}^{n})^{2}}\right) 
.$\ \par
\ \par
It follows that:\ \par
$v_{i+1}^{n+\frac{1}{2}}-v_{i}^{n+\frac{1}{2}}=(v_{i+1}^{n}-v_{i}^{n})(1-\mu 
_{i}^{n})+\mu _{i-1}^{n}(v_{i}^{n}-v_{i-1}^{n})$\ \par
\ \par
$-2\lambda \Delta t_{n}x_{i}(v_{i+1}^{n}+v_{i}^{n})(v_{i+1}^{n}-v_{i}^{n})-\lambda 
\Delta t_{n}h_{i+1}(v_{i+1}^{n})^{2}-\lambda \Delta t_{n}h_{i}(v_{i}^{n})^{2}$\ 
\par
\ \par
with  $\displaystyle \mu _{i}^{n}=\lambda \frac{\displaystyle \Delta t_{n}}{\displaystyle h_{i+1}}x_{i}^{2}(v_{i+1}^{n}+v_{i}^{n}),\ 
1\leq i\leq N-1,\ \mu _{N}^{n}=0.$ \ \par
\ \par
From the stability condition ~(\ref{MarchesIncomplets28}), we get $0\leq \mu _{i}^{n}\leq 1$ 
and we deduce \ \par

\begin{displaymath} 
Var\ (v_{h}^{n+\frac{1}{2}};{\mathbb{R}}^{+})\leq Var(v_{h}^{n};{\mathbb{R}}^{+})\ 
(1+4\lambda \Delta t_{n}\left\Vert{\ \hat v_{h}^{n}}\right\Vert _{L^{\infty 
}({\mathbb{R}}^{+})})+2\lambda \Delta t_{n}\left\Vert{v_{h}^{n}}\right\Vert 
_{L^{\infty }({\mathbb{R}}^{+})}\left\Vert{v_{h}^{n}}\right\Vert _{L^{1}({\mathbb{R}}^{+})}.\label{MarchesIncomplets38}\end{displaymath} 
\ \par
We have an analogous estimate for $\lambda <0$.\ \par
\ \par
Furthermore, from ~(\ref{MarchesIncomplets21}), we obtain the following 
equalities:\ \par

\begin{displaymath} 
v_{i+1}^{n+1}-v_{i}^{n+1}+\Delta t_{n}\left({\frac{\alpha _{i}+(1-\delta 
_{i})\gamma _{i}}{h_{i}}+\frac{\alpha _{i}-\delta _{i}\gamma _{i}}{h_{i+1}}-\frac{\gamma 
_{i+1}-\gamma _{i}}{h_{i+1}}}\right) \left({v_{i+1}^{n+1}-v_{i}^{n+1}}\right) 
\end{displaymath} \ \par
\ \par

\begin{displaymath} 
-\Delta t_{n}\frac{\alpha _{i+1}+(1-\delta _{i+1})\gamma _{i+1}}{h_{i+1}}\left({v_{i+2}^{n+1}-v_{i+1}^{n+1}}\right) 
-\Delta t_{n}\frac{\alpha _{i-1}-\delta _{i-1}\gamma _{i-1}}{h_{i}}\left({v_{i}^{n+1}-v_{i-1}^{n+1}}\right) 
\end{displaymath} \ \par

\begin{displaymath} 
=v_{i+1}^{n+\frac{1}{2}}-v_{i}^{n+\frac{1}{2}}+\Delta t_{n}\left({\frac{\gamma 
_{i+1}-\gamma _{i}}{h_{i+1}}-\frac{\gamma _{i}-\gamma _{i-1}}{h_{i}}}\right) 
v_{i}^{n+1}-\Delta t_{n}\left({\frac{F_{1},_{i+1}-F_{1i}}{h_{i+1}}-\frac{F_{1i}-F1,_{i-1}}{h_{i}}}\right) 
,\end{displaymath} \ \par
for $1\leq i\leq N-1$.\ \par

\begin{displaymath} 
v_{N+1}^{n+1}-v_{N}^{n+1}+\Delta t_{n}\left({\frac{\alpha _{N}}{h_{N}}(v_{N+1}^{n+1}-v_{N}^{n+1})-\frac{\alpha 
_{N-1}-\gamma _{N-1}}{h_{N}}(v_{N}^{n+1}-v_{N-1}^{n+1})}\right) \end{displaymath} 
\ \par

\begin{displaymath} 
=v_{N+1}^{n+\frac{1}{2}}-v_{N}^{n+\frac{1}{2}}-\Delta t_{n}\frac{\gamma 
_{N}-\gamma _{N-1}}{h_{N}}+\Delta t_{n}\frac{F_{1N}-F_{1N-1}}{h_{N}}.\end{displaymath} 
\ \par
It follows that:\ \par

\begin{displaymath} 
\sum_{i=1}^{N-1}{}\left({1-\Delta t_{n}\frac{\gamma _{i+1}-\gamma _{i}}{h_{i+1}}}\right) 
\left\vert{v_{i+1}^{n+1}-v_{i}^{n+1}}\right\vert +v_{N}^{n+1}\leq Var(v_{h}^{n+\frac{1}{2}};{\mathbb{R}}^{+})\end{displaymath} 
\ \par

\begin{displaymath} 
+\Delta t_{n}\displaystyle \sum_{i=1}^{N-1}{\displaystyle \left\vert{\displaystyle 
\frac{\displaystyle \gamma _{i+1}-\gamma _{i}}{\displaystyle h_{i+1}}-\frac{\displaystyle 
\gamma _{i}-\gamma _{i-1}}{\displaystyle h_{i}}}\right\vert v_{i}^{n+1}}\end{displaymath} 
\ \par

\begin{displaymath} 
+\Delta t_{n}\displaystyle \sum_{i=1}^{N-1}{\displaystyle \left\vert{\displaystyle 
\frac{\displaystyle F_{1i}-F_{1},_{i-1}}{\displaystyle h_{i+1}}-\frac{\displaystyle 
F_{1i}-F_{1},_{i-1}}{\displaystyle h_{i}}}\right\vert }+\Delta t_{n}\displaystyle 
\left\vert{\displaystyle \frac{\displaystyle F_{1N}-F_{1},_{N-1}}{\displaystyle 
h_{N}}}\right\vert \ .\end{displaymath} \ \par
\ \par
It results from ~(\ref{MarchesIncomplets33}): $\displaystyle 1-\Delta t_{n}\frac{\displaystyle \gamma _{i+1}-\gamma _{i}}{\displaystyle 
h_{i+1}}\geq 1-c_{1}\Delta t_{n}$.\ \par
Besides, we have:  $\displaystyle \displaystyle \sum_{i=1}^{N-1}{\displaystyle \left\vert{\displaystyle 
\frac{\displaystyle \gamma _{i+1}-\gamma _{i}}{\displaystyle h_{i+1}}-\frac{\displaystyle 
\gamma _{i}-\gamma _{i-1}}{\displaystyle h_{i}}}\right\vert v_{i}^{n+1}}\leq 
Var(g_{1}';{\mathbb{R}}^{+})\displaystyle \left\Vert{\displaystyle v_{h}^{n+1}}\right\Vert 
_{\displaystyle L^{\infty }({\mathbb{R}}^{+})}$\ \par
\ \par
and $\displaystyle \displaystyle \sum_{i=1}^{N-1}{\displaystyle \left\vert{\displaystyle 
\frac{\displaystyle F_{1,i+1}-F_{1i}}{\displaystyle h_{i+1}}-\frac{\displaystyle 
F_{1i}-F_{1,i-1}}{\displaystyle h_{i}}}\right\vert }\leq Var(F_{1}';{\mathbb{R}}^{+})$. 
\ \par
\ \par
It follows that:
\begin{displaymath} 
(1-c_{1}\Delta t_{n})\ Var(v_{h}^{n+1};{\mathbb{R}}^{+})\leq Var(v_{h}^{n+\frac{1}{2}};{\mathbb{R}}^{+})+C\Delta 
t_{n}\left\Vert{v_{h}^{n+1}}\right\Vert _{L^{1}({\mathbb{R}}^{+})}+\Delta 
t_{n}Var(F_{1}',{\mathbb{R}}^{+})\end{displaymath} \ \par
and with ~(\ref{MarchesIncomplets38}), we get:\ \par

\begin{displaymath} 
(1-c_{1}\Delta t_{n})Var(v_{h}^{n+1};{\mathbb{R}}^{+}{\mathrm{)}}\leq 
Var(v_{h}^{n};{\mathbb{R}}^{+})(1+4\lambda \Delta t_{n}\displaystyle 
\left\Vert{\displaystyle \hat v_{h}^{n}}\right\Vert _{\displaystyle L^{\infty 
}({\mathbb{R}}^{+})})\end{displaymath} \ \par
$+2\lambda \Delta t_{n}\left\Vert{v_{h}^{n}}\right\Vert _{L^{\infty }({\mathbb{R}}^{+})}\left\Vert{v_{h}^{n}}\right\Vert 
_{L^{1}({\mathbb{R}}^{+})}+C\Delta t_{n}\left\Vert{v_{h}^{n+1}}\right\Vert 
_{L^{1}({\mathbb{R}}^{+})}$ $+\Delta t_{n}Var(F_{1}';{\mathbb{R}}^{+}).$\ 
\par
\ \par
This concludes the proof.\ \par
\ \par
\subsection{Convergence of the scheme and uniqueness of the solution}
\ \par
From all these estimates, we can deduce the convergence of the numerical 
solution to a weak solution and we prove that this solution is unique.\ 
\par
\ \par
{\hskip 1.8em}A function $u_{1}$ is called a {\it{weak solution}} of 
~(\ref{MarchesIncomplets6}) if $u_{1}\in C([0,T]);\ W^{1,\infty }({\mathbb{R}}^{+})),$\ 
\par
$\ x\frac{\partial u_{1}}{\partial x}\in C(0,T;L^{\infty }({\mathbb{R}}+))$, 
and\ \par

\begin{displaymath} 
\int_{0}^{T}{\int_{{\mathbb{R}}^{+}}^{}{u_{1}\frac{\partial \phi }{\partial 
t}dx}dt}-\int_{{\mathbb{R}}^{+}}^{}{u_{1}(T)\phi (T)dx}-\frac{k^{2}}{2}\int_{0}^{T}{\int_{{\mathbb{R}}^{+}}^{}{\frac{\partial 
u_{1}}{\partial x}\frac{\partial }{\partial x}(x^{2}\phi )dx}dx}\end{displaymath} 
 \ \par
\ \par

\begin{displaymath} 
+\displaystyle \int_{0}^{T}{\displaystyle \int_{{\mathbb{R}}^{+}}^{}{g_{1}\frac{\displaystyle 
\partial u_{1}}{\displaystyle \partial x}\phi dx}dt}-\lambda \displaystyle 
\int_{0}^{T}{\displaystyle \int_{{\mathbb{R}}^{+}}^{}{x^{2}\displaystyle 
\left({\displaystyle \frac{\displaystyle \partial u_{1}}{\displaystyle 
\partial x}}\right) ^{2}\phi dx}dt}=\displaystyle \int_{0}^{T}{\displaystyle 
\int_{{\mathbb{R}}^{+}}^{}{F_{1}\phi dx}dt}\ \end{displaymath} \ \par
for any function $\phi $ with a compact support in $[0,T]\times {\mathbb{R}}^{+},\phi \in C^{1}(0,T\times {\mathbb{R}}^{+}).$\ 
\par
\ \par
\begin{Theorem} Problem ~(\ref{MarchesIncomplets6}) admits at most one 
weak solution.\ \par
\end{Theorem}
\ \par
{\it{Proof: }}Let $u_{1}$ and $\hat u_{1}$ two weak solutions of ~(\ref{MarchesIncomplets6}). 
We denote $w=u_{1}-\hat u_{1}$.The function $w$ satisfies: 
\begin{equation} 
\frac{\displaystyle \partial w}{\displaystyle \partial t}-\frac{\displaystyle 
1}{\displaystyle 2}k^{2}x^{2}\frac{\displaystyle \partial ^{2}w}{\displaystyle 
\partial x^{2}}-g_{1}(x)\frac{\displaystyle \partial w}{\displaystyle 
\partial x}+\lambda x^{2}\frac{\displaystyle \partial w}{\displaystyle 
\partial x}\displaystyle \left({\displaystyle \frac{\displaystyle \partial 
u_{1}}{\displaystyle \partial x}+\frac{\displaystyle \partial \hat u_{1}}{\displaystyle 
\partial x}}\right) =0.\label{MarchesIncomplets41}
\end{equation} \ \par
Let us denote by $\psi $ a function in $C^{1}({\mathbb{R}}^{+})$ satisfying:\ 
\par
$\displaystyle 0\leq \psi (x)\leq 1,\ \psi (x)=1\ if\ 0\leq x\leq 1,\ \psi (x)=0\ if\ 
x\geq 2$ and $\psi \ $decreasing on $(1,2)$ \ \par
and we define $\psi _{\nu }(x)=\psi (\frac{x}{\nu }),\ x\geq 0,\ \nu >0$.\ 
\par
By multiplying ~(\ref{MarchesIncomplets41}) by $\displaystyle \psi _{\nu }\frac{\displaystyle w}{\displaystyle (1+x)^{2}}$, 
and integrating on ${\mathbb{R}}^{+}$, we get:\ \par

\begin{displaymath} 
\frac{\displaystyle 1}{\displaystyle 2}\frac{\displaystyle d}{\displaystyle 
dt}\displaystyle \left({\displaystyle \displaystyle \int_{0}^{+\infty 
}{w^{2}\frac{\displaystyle \psi _{\nu }}{\displaystyle (1+x)^{2}}dx}}\right) 
+\frac{\displaystyle 1}{\displaystyle 2}k^{2}\displaystyle \int_{0}^{+\infty 
}{\displaystyle \left({\displaystyle \frac{\displaystyle \partial w}{\displaystyle 
\partial x}}\right) ^{2}\frac{\displaystyle x^{2}}{\displaystyle (1+x)^{2}}\psi 
_{\nu }dx}+\frac{\displaystyle 1}{\displaystyle 2}\displaystyle \int_{0}^{+\infty 
}{w^{2}\frac{\displaystyle g_{1}(x)}{\displaystyle (1+x)^{2}}\frac{\displaystyle 
d\psi _{\nu }}{\displaystyle dx}dx}\end{displaymath} \ \par

\begin{displaymath} 
=-\frac{\displaystyle k^{2}}{\displaystyle 2}\displaystyle \int_{0}^{+\infty 
}{\frac{\displaystyle \partial w}{\displaystyle \partial x}w\frac{\displaystyle 
x^{2}}{\displaystyle (1+x)^{2}}\frac{\displaystyle d\psi _{\nu }}{\displaystyle 
dx}dx}-k^{2}\displaystyle \int_{0}^{+\infty }{\frac{\displaystyle \partial 
w}{\displaystyle \partial x}w\frac{\displaystyle x}{\displaystyle (1+x)^{3}}\psi 
_{\nu }dx}\end{displaymath} \ \par

\begin{equation} 
-\frac{\displaystyle 1}{\displaystyle 2}\displaystyle \int_{0}^{+\infty 
}{w^{2}\psi _{\nu }}\frac{\displaystyle d}{\displaystyle dx}\displaystyle 
\left({\displaystyle \frac{\displaystyle g_{1}(x)}{\displaystyle (1+x)^{2}}}\right) 
dx-\lambda \displaystyle \int_{0}^{+\infty }{\frac{\displaystyle \partial 
w}{\displaystyle \partial x}w\displaystyle \left({\displaystyle \frac{\displaystyle 
\partial u_{1}}{\displaystyle \partial x}+\frac{\displaystyle \partial 
\hat u_{1}}{\displaystyle \partial x}}\right) \frac{\displaystyle x^{2}}{\displaystyle 
(1+x)^{2}}\psi _{\nu }dx}{\mathrm{.}}\label{MarchesIncomplets42}
\end{equation} \ \par
\ \par
We estimate now each term of this equality.\ \par
We have: $\displaystyle \displaystyle \int_{0}^{+\infty }{w^{2}\frac{\displaystyle g_{1}(x)}{\displaystyle 
(1+x)^{2}}\frac{\displaystyle d\psi _{\nu }}{\displaystyle dx}dx}\geq 
0$ if $\nu \geq \sigma _{2}$ since $g_{1}(x)\leq 0\ $and $\psi '_{\nu }(x)\leq 0$.\ 
\par
\ \par
Besides, we have: $\displaystyle \psi '_{\nu }(x)=\frac{\displaystyle 1}{\displaystyle \nu }\psi '(\frac{\displaystyle 
x}{\displaystyle \nu })$ and since $x\frac{\partial w}{\partial x}\in C(0,T;L^{\infty }({\mathbb{R}}^{+})),\ 
w\in C(0,T;L^{\infty }({\mathbb{R}}^{+}))$, we obtain:\ \par
\ \par
$\displaystyle \frac{\displaystyle k^{2}}{\displaystyle 2}\displaystyle \left\vert{\displaystyle 
\displaystyle \int_{0}^{+\infty }{\frac{\displaystyle \partial w}{\displaystyle 
\partial x}w\frac{\displaystyle x^{2}}{\displaystyle (1+x)^{2}}\frac{\displaystyle 
d\psi _{\nu }}{\displaystyle dx}dx}}\right\vert \leq \frac{\displaystyle 
C}{\displaystyle \nu }$ ($C$ depending on $\left\Vert{x\frac{\partial w}{\partial x}}\right\Vert _{L^{\infty }({\mathbb{R}}^{+})}\ 
$and $\left\Vert{w}\right\Vert _{L^{\infty }({\mathbb{R}}^{+})}$).\ \par
\ \par
\ \par
We estimate the second term of the second member of ~(\ref{MarchesIncomplets42}) 
and we get:\ \par

\begin{displaymath} 
k^{2}\left\vert{\int_{0}^{+\infty }{\frac{\partial w}{\partial x}w\frac{x}{(1+x)^{3}}\psi 
_{\nu }dx}}\right\vert \leq k^{2}\left({\int_{0}^{+\infty }{\frac{w^{2}}{(1+x)^{2}}\psi 
_{\nu }dx}}\right) ^{\frac{1}{2}}\left({\int_{0}^{+\infty }{\left({\frac{\partial 
w}{\partial x}}\right) ^{2}\frac{x^{2}}{(1+x)^{4}}\psi _{\nu }dx}}\right) 
^{\frac{1}{2}}\end{displaymath} \ \par

\begin{displaymath} 
\leq \alpha k^{2}\int_{0}^{+\infty }{\left({\frac{\partial w}{\partial 
x}}\right) ^{2}\frac{x^{2}}{(1+x)^{2}}\psi _{\nu }dx}+\frac{k^{2}}{4\alpha 
}\int_{0}^{+\infty }{w^{2}\frac{\psi _{\nu }}{(1+x)^{2}}dx},\end{displaymath} 
\ \par
$\alpha >0$.\ \par
\ \par
Further, we have:   $\displaystyle \displaystyle \left\vert{\displaystyle \displaystyle \int_{0}^{+\infty 
}{w^{2}\psi _{\nu }\frac{\displaystyle d}{\displaystyle dx}\displaystyle 
\left({\displaystyle \frac{\displaystyle g_{1}(x)}{\displaystyle (1+x)^{2}}}\right) 
dx}}\right\vert \leq \displaystyle \int_{0}^{+\infty }{w^{2}\frac{\displaystyle 
\psi _{\nu }}{\displaystyle (1+x)^{2}}\displaystyle \left\vert{\displaystyle 
\frac{\displaystyle g_{1}'(x)(1+x)-2g_{1}(x)}{\displaystyle 1+x}}\right\vert 
dx}$. \ \par
\ \par
From the hypotheses on the function $g_{1}$, we get\ \par
\ \par
$\displaystyle \displaystyle \left\vert{\displaystyle \frac{\displaystyle g_{1}'(x)(1+x)-2g_{1}(x)}{\displaystyle 
1+x}}\right\vert \leq \displaystyle \left\vert{\displaystyle \frac{\displaystyle 
\phi (x)+x\phi '(x)}{\displaystyle 1+x}}\right\vert +\displaystyle \left\vert{\displaystyle 
\phi (x)-x\phi '(x)}\right\vert $\ \par
and this quantity is  bounded.Then , we obtain:\ \par
$\displaystyle \displaystyle \left\vert{\displaystyle \displaystyle \int_{0}^{+\infty 
}{w^{2}\psi _{\nu }\frac{\displaystyle d}{\displaystyle dx}\displaystyle 
\left({\displaystyle \frac{\displaystyle g_{1}(x)}{\displaystyle (1+x)^{2}}}\right) 
dx}}\right\vert \leq C\displaystyle \int_{0}^{+\infty }{w^{2}\frac{\displaystyle 
\psi _{\nu }}{\displaystyle (1+x)^{2}}dx}.$\ \par
\ \par
It remains to study the last term of ~(\ref{MarchesIncomplets42}):\ \par
Since $\displaystyle x\frac{\displaystyle \partial u_{1}}{\displaystyle \partial x}$ 
and $\displaystyle x\frac{\displaystyle \partial \hat u_{1}}{\displaystyle \partial x}\in 
C(0,T;L^{\infty }({\mathbb{R}}^{+}))$, we get:\ \par
\ \par

\begin{displaymath} 
\left\vert{\int_{0}^{+\infty }{\frac{\partial w}{\partial x}w\left({\frac{\partial 
u_{1}}{\partial x}+\frac{\partial \hat u_{1}}{\partial x}}\right) \frac{x^{2}}{(1+x)^{2}}\psi 
_{\nu }dx}}\right\vert \leq C\int_{0}^{+\infty }{\left\vert{\frac{\partial 
w}{\partial x}}\right\vert \left\vert{w}\right\vert \frac{x}{(1+x)^{2}}\psi 
_{\nu }dx}\end{displaymath} \ \par

\begin{displaymath} 
\leq \alpha k^{2}\int_{0}^{+\infty }{\left({\frac{\partial w}{\partial 
x}}\right) ^{2}\frac{x^{2}}{(1+x)^{2}}\psi _{\nu }dx}+\frac{C^{2}}{4\alpha 
k^{2}}\int_{0}^{+\infty }{w^{2}\frac{\psi _{\nu }}{(1+x)^{2}}dx}.\end{displaymath} 
\ \par
\ \par
We deduce from all these estimates:\ \par

\begin{displaymath} 
\frac{1}{2}\frac{d}{dt}\left({\int_{0}^{+\infty }{w^{2}\frac{\psi _{\nu 
}}{(1+x)^{2}}dx}}\right) +k^{2}(\frac{1}{2}-2\alpha )\int_{0}^{+\infty 
}{\left({\frac{\partial w}{\partial x}}\right) ^{2}\frac{x^{2}}{(1+x)^{2}}\psi 
_{\nu }dx}\end{displaymath} \ \par

\begin{displaymath} 
\leq \frac{C}{\nu }+C_{1}\int_{0}^{+\infty }{w^{2}\frac{\psi _{\nu }}{(1+x)^{2}}dx}.\end{displaymath} 
\ \par
\ \par
If we choose $\alpha <\frac{1}{4}$, we obtain by using the  Gronwall's 
lemma:\ \par
$\displaystyle \displaystyle \int_{0}^{+\infty }{w^{2}\frac{\displaystyle \psi _{\nu 
}}{\displaystyle (1+x)^{2}}dx}\leq \frac{\displaystyle 2CT}{\displaystyle 
\nu }e^{\displaystyle 2C_{1}t}$, $0\leq t\leq T$\ \par
and if $\nu {\longrightarrow}+\infty $, we get $\displaystyle \displaystyle \int_{0}^{+\infty }{\frac{\displaystyle w^{2}}{\displaystyle 
(1+x)^{2}}dx=0}$\ \par
and we deduce $w=0$ and the problem admits at most one solution.\ \par
\ \par
\ \par
We prove now the convergence of the numerical solution to this weak solution 
and thus, we obtain the existence of a solution.\ \par
We define the functions $u_{1h\Delta t}$ and $v_{1h\Delta t}$ by\ \par
$\displaystyle u_{1h\Delta t}=u_{h}^{n}+\frac{\displaystyle t-t_{n}}{\displaystyle \Delta 
t_{n}}(u_{1h}^{n+1}-u_{1h}^{n}),$\ \par
\ \par
$\displaystyle v_{h\Delta t}=v_{h}^{n}+\frac{\displaystyle t-t_{n}}{\displaystyle \Delta 
t_{n}}(v_{h}^{n+1}-v_{h}^{n})$, $t_{n}\leq t\leq t_{n+1}$.\ \par
\ \par
\begin{Theorem} Assume that the hypotheses of proposition ~\ref{MarchesIncomplets36} 
are satisfied. The sequence $u_{1h\Delta t}$ converges uniformly to the 
weak solution of ~(\ref{MarchesIncomplets6}) on any compact of $[0,T]\times {\mathbb{R}}^{+}$.\ 
\par
\end{Theorem}
\ \par
{\it{Proof:}} The functions $(u_{1h\Delta t})\ $ are uniformly bounded 
in $C(0,T;W^{1,\infty }({\mathbb{R}}^{+}))$.\ \par
Further, for $R>0$, we get the estimate:\ \par
\ \par
$\displaystyle \displaystyle \left\Vert{\displaystyle \frac{\displaystyle u_{1h}^{n+1}-u_{1h}^{n}}{\displaystyle 
\Delta t_{n}}}\right\Vert _{\displaystyle L^{1}(0,R)}\leq \frac{\displaystyle 
k^{2}R^{2}}{\displaystyle 2}Var(v_{h}^{n+1};{\mathbb{R}}^{+})+C(g_{1},R)\displaystyle 
\left\Vert{\displaystyle v_{h}^{n+1}}\right\Vert _{\displaystyle L^{1}({\mathbb{R}}^{+})}$\ 
\par
$+\left\vert{\lambda }\right\vert R\left\Vert{\hat v_{h}^{n}}\right\Vert 
_{L^{\infty }({\mathbb{R}}^{+})}\left\Vert{v_{h}^{n}}\right\Vert _{L^{1}({\mathbb{R}}^{+})}+\left\Vert{F_{1}}\right\Vert 
_{L^{1}({\mathbb{R}}^{+})}$\ \par
\ \par
with $\displaystyle C(g_{1},R)=\underset{\displaystyle x\leq R}{\displaystyle {\mathrm{sup}}}\displaystyle 
\left\vert{\displaystyle g_{1}(x)}\right\vert $. \ \par
Thus, the time derivatives of $u_{h\Delta t}$ are uniformly bounded in 
$L^{\infty }(0,T;L^{1}(0,R))$ for any $R>0$.\ \par
So, we can extract from the sequence $(u_{h\Delta t})$ a subsequence, 
again labeled $u_{h\Delta t}$ which converges uniformly  on any compact 
subset of $[0,T]\times {\mathbb{R}}^{+}$ to a function $u_{1}$ ~\cite{Simon}.\ 
\par
{\hskip 1.8em}The functions $v_{h\Delta t}$ are uniformly bounded in 
$C(0,T;BV({\mathbb{R}}^{+}))$;\ \par
besides, we have:\ \par

\begin{equation} 
\displaystyle \left\Vert{\displaystyle v_{h}^{n+1}-v_{h}^{n}}\right\Vert 
_{\displaystyle H^{-2}(0,R)}=\underset{\displaystyle \phi \in H_{0}^{2}(0,R)}{\displaystyle 
{\mathrm{sup}}}\frac{\displaystyle <v_{h}^{n+1}-v_{h}^{n},\phi >}{\displaystyle 
\displaystyle \left\Vert{\displaystyle \phi }\right\Vert _{\displaystyle 
H_{0}^{2}(0,R)}}=\ \underset{\displaystyle \phi \in H_{0}^{2}(0,R)}{\displaystyle 
{\mathrm{sup}}}\frac{\displaystyle <u_{h}^{n+1}-u_{h}^{n},\phi _{x}>}{\displaystyle 
\displaystyle \left\Vert{\displaystyle \phi }\right\Vert _{\displaystyle 
H_{0}^{2}(0,R)}}
\end{equation} \ \par
and $\displaystyle \displaystyle \left\Vert{\displaystyle \frac{\displaystyle v_{h}^{n+1}-v_{h}^{n}}{\displaystyle 
\Delta t_{n}}}\right\Vert _{\displaystyle H^{-2}(0,R)}\leq C\displaystyle 
\left\Vert{\displaystyle \frac{\displaystyle u_{h}^{n+1}-u_{h}^{n}}{\displaystyle 
\Delta t_{n}}}\right\Vert _{\displaystyle L^{1}(0,R)}\leq C(R).$\ \par
Then, we can extract from $(v_{h\Delta t})$ a subsequence, again labeled 
$v_{h\Delta t}$ which converges to a function $v=\frac{\partial u_{1}}{\partial x}$ 
in $C(0,T;L^{1}(Q))$ for any compact $Q\subset {\mathbb{R}}^{+}$ ~\cite{Simon}. 
\ \par
\ \par
We get easily  that $u_{1}$ is a weak solution. Since this solution is 
unique, all the sequence is converging to $u_{1}$.\ \par
\ \par
So, we have obtained the following result:\ \par
\begin{Theorem} Problem ~(\ref{MarchesIncomplets6}) admits a unique weak 
solution.\ \par
\end{Theorem}
\ \par

 \ \par
The following figure represents $a(0)$ with a time maturity~\cite{Bolley} 
equal to 1 for different values of $\rho \cdot $ the values of the other 
parameters are those proposed in ~\cite{Kogan}: $k=0.4,\ \delta =2,\ \sigma _{1}=0.153,\ \mu =0.7,\ \sigma _{0}=0.01$.\ 
\par

\begin{figure}[h]
\begin{center}
\rotatebox{0}{\resizebox{10cm}{!}{\includegraphics{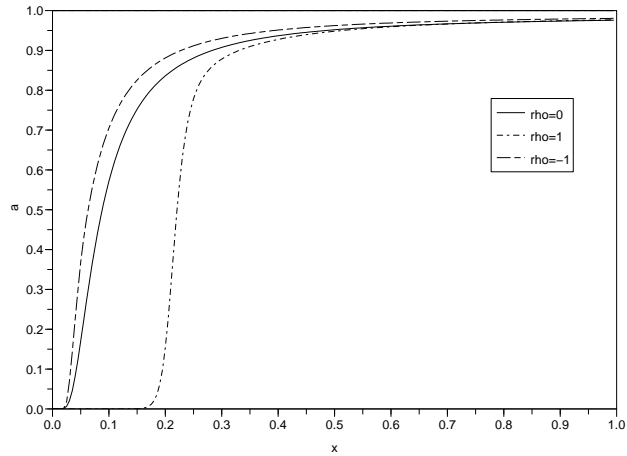}}}
\caption{a with T=1}
\label{MarchesIncomplets100}
\end{center}
\end{figure}
\ \par

\section{Computation of $u_{2}$}
\setcounter{equation}{0}$u_{2}$ is solution of \ \par
\ \par

\begin{displaymath} 
\frac{\partial u_{2}}{\partial t}-\frac{k^{2}x^{2}}{2}\frac{\partial ^{2}u_{2}}{\partial 
x^{2}}-\frac{x^{2}y^{2}}{2}\frac{\partial ^{2}u_{2}}{\partial y^{2}}-\rho 
kx^{2}y\frac{\partial ^{2}u_{2}}{\partial x\partial y}-g_{2}(x)\frac{\partial 
u_{2}}{\partial x}\end{displaymath} \ \par
\ \par

\begin{displaymath} 
-(1-\rho ^{2})k^{2}x^{2}\frac{\partial u_{1}}{\partial x}\frac{\partial 
u_{2}}{\partial x}=0.\end{displaymath} \ \par
The initial condition is the payoff function. If the European option 
is a put, the payoff function is given by:$\displaystyle F(y)={\mathrm{Max}}(E-y,0)$ 
if $E$ is the exercise price.\ \par
\ \par
\subsection{Existence and uniqueness of the solution}
\ \par
In order to obtain a bilinear form satisfying Garding's inequality, we 
make a change of unknown.\ \par
\ \par
We denote:$\displaystyle \ \hat u_{2}=e^{-\alpha x}u_{2}$, $\alpha >0$.\ 
\par
The function $\hat u_{2}$ is solution of:\ \par
\ \par

\begin{displaymath} 
\frac{\displaystyle \partial \hat u_{2}}{\displaystyle \partial t}-\frac{\displaystyle 
k^{2}}{\displaystyle 2}x^{2}\frac{\displaystyle \partial ^{2}\hat u_{2}}{\displaystyle 
\partial x^{2}}-\frac{\displaystyle x^{2}y^{2}}{\displaystyle 2}\frac{\displaystyle 
\partial ^{2}\hat u_{2}}{\displaystyle \partial y^{2}}-\rho kx^{2}y\frac{\displaystyle 
\partial ^{2}\hat u_{2}}{\displaystyle \partial x\partial y}-\frac{\displaystyle 
\partial \hat u_{2}}{\displaystyle \partial x}\displaystyle \left({\displaystyle 
\alpha k^{2}x^{2}+g_{2}(x)+(1-\rho ^{2})k^{2}x^{2}\frac{\displaystyle 
\partial u_{1}}{\displaystyle \partial x}}\right) \end{displaymath}  
$\displaystyle -\rho \alpha kx^{2}y\frac{\displaystyle \partial \hat u_{2}}{\displaystyle 
\partial y}-\displaystyle \left({\displaystyle \frac{\displaystyle \alpha 
^{2}k^{2}}{\displaystyle 2}x^{2}+\alpha g_{2}(x)+(1-\rho ^{2})\alpha 
k^{2}x^{2}\frac{\displaystyle \partial u_{1}}{\displaystyle \partial 
x}}\right) \hat u_{2}=0,\label{1MarchesIncomplets0}$\ \par
\ \par
with the initial condition: $\displaystyle \hat u_{2}(0)=e^{-\alpha x}F(y)$.\ 
\par
\ \par
We define a variational formulation of this problem.\ \par
\ \par
Let us consider the following space $\hat V_{2}$ defined by:\ \par
$\displaystyle \hat V_{2}=\displaystyle \left\lbrace{\displaystyle \displaystyle \left.{\displaystyle 
v\in {\mathcal{D}}'(\Omega )/\ v\in L^{2}(\Omega ),\ xv\in L^{2}(\Omega 
),\ xv_{x}\in L^{2}(\Omega ),\ xyv_{y}\in L^{2}(\Omega )}\right\rbrace 
}\right. $\ \par
with $\Omega ={\mathbb{R}}^{+}\times {\mathbb{R}}^{+}$.\ \par
\ \par
This space with the norm: \ \par
$\displaystyle \displaystyle \left\Vert{\displaystyle v}\right\Vert =\displaystyle \left({\displaystyle 
\displaystyle \int_{\Omega }^{}{v^{2}+x^{2}v^{2}+x^{2}\displaystyle \left({\displaystyle 
\frac{\displaystyle \partial v}{\displaystyle \partial x}}\right) ^{2}+x^{2}y^{2}\displaystyle 
\left({\displaystyle \frac{\displaystyle \partial v}{\displaystyle \partial 
y}}\right) ^{2}}}\right) ^{\displaystyle \frac{\displaystyle 1}{\displaystyle 
2}}$\ \par
is a Hilbert space and ${\mathcal{D}}(\Omega )$ is dense in$\hat V_{2}$ 
~\cite{Bolley} .\ \par
Besides, we have the estimates: \ \par
$\displaystyle \displaystyle \left\Vert{\displaystyle xv}\right\Vert _{\displaystyle 
L^{2}(\Omega )}\leq 2\displaystyle \left\Vert{\displaystyle xy\frac{\displaystyle 
\partial v}{\displaystyle \partial y}}\right\Vert _{\displaystyle L^{2}(\Omega 
)}$,  $\displaystyle \displaystyle \left\Vert{\displaystyle v}\right\Vert _{\displaystyle L^{2}(\Omega 
)}\leq 2\displaystyle \left\Vert{\displaystyle x\frac{\displaystyle \partial 
v}{\displaystyle \partial x}}\right\Vert _{\displaystyle L^{2}(\Omega 
)}$\ \par
and the semi-norm:\ \par

\begin{displaymath} 
\left\Vert{v}\right\Vert _{2}=\ \left({\left\Vert{x\frac{\partial v}{\partial 
x}}\right\Vert ^{2}_{L^{2}(\Omega )}+\left\Vert{xy\frac{\partial v}{\partial 
y}}\right\Vert ^{2}_{L^{2}(\Omega )}}\right) ^{2}\end{displaymath} \ 
\par
is a norm in $\hat V_{2}$ equivalent to the norm $\left\Vert{.}\right\Vert _{\hat V_{2}}$ 
~\cite{Bolley} .\ \par
\ \par
We define on $\hat V_{2}\times \hat V_{2}$ the bilinear form $\hat b$ 
by : \ \par
\ \par
$\displaystyle \forall u,v\in \hat V_{2},\ \hat b(u,v)=\frac{\displaystyle k^{2}}{\displaystyle 
2}\displaystyle \int_{\Omega }^{}{\frac{\displaystyle \partial u}{\displaystyle 
\partial x}\frac{\displaystyle \partial }{\displaystyle \partial x}(x^{2}v)dxdy}+\frac{\displaystyle 
1}{\displaystyle 2}\displaystyle \int_{\Omega }^{}{x^{2}\frac{\displaystyle 
\partial u}{\displaystyle \partial y}\frac{\displaystyle \partial }{\displaystyle 
\partial y}(y^{2}v)dxdy}$\ \par
$\displaystyle +\rho k\displaystyle \int_{\Omega }^{}{x^{2}\frac{\displaystyle \partial 
u}{\displaystyle \partial x}\frac{\displaystyle \partial }{\displaystyle 
\partial y}(yv)dxdy}-\displaystyle \int_{\Omega }^{}{(\alpha k^{2}x^{2}+g_{2}(x)+(1-\rho 
^{2})k^{2}x^{2}\frac{\displaystyle \partial u_{1}}{\displaystyle \partial 
x})\frac{\displaystyle \partial u}{\displaystyle \partial x}vdxdy\ }$\ 
\par
$\displaystyle -\rho \alpha k\displaystyle \int_{\Omega }^{}{x^{2}y\frac{\displaystyle 
\partial u}{\displaystyle \partial y}vdxdy}\ -\displaystyle \int_{\Omega 
}^{}{\displaystyle \left({\displaystyle \frac{\displaystyle \alpha ^{2}k^{2}}{\displaystyle 
2}x^{2}+\alpha g_{2}(x)+(1-\rho ^{2})\alpha k^{2}x^{2}\frac{\displaystyle 
\partial u_{1}}{\displaystyle \partial x}}\right) uvdxdy}$.\ \par
\ \par
\ \par
We have proved in the preceding section that the function $\displaystyle 
x\frac{\displaystyle \partial u_{1}}{\displaystyle \partial x}\ $is in 
$C(0,T;L^{\infty }({\mathbb{R}}^{+}))$.\ \par
Besides,  we have $g_{2}(x)=-\delta x(x-\sigma _{1})-\rho kxf(x)$.\ \par
We denote $f_{1}(x)=xf(x)$ and we assume that $f_{1},f_{1}'\in L^{\infty }({\mathbb{R}}^{+})$.\ 
\par
It is clear that the the bilinear form $\hat b$ is continue on $\hat V_{2}\times \hat V_{2}$.\ 
\par
\ \par
We prove now that it is possible to choose $\alpha >0$ such that this 
bilinear form satisfies Garding's inequality.\ \par
\ \par
\begin{Proposition} \label{1MarchesIncomplets1}If $\left\vert{\rho }\right\vert <1$ 
and $k<\frac{2}{3}\delta $, there exists positive constants $C,\ c,\ \alpha $ 
depending on $\delta ,\sigma _{1},k,\rho $ such that 
\begin{displaymath} 
\forall v\in \hat V_{2},\ \hat b(v,v)\geq C\left\Vert{v}\right\Vert _{2}^{2}-c\left\Vert{v}\right\Vert 
^{2}_{L^{2}(\Omega )}.\end{displaymath} \ \par
\end{Proposition}
{\it{Proof: }}We have the equality:\ \par
$\displaystyle \hat b(v,v)=\displaystyle \sum_{i=1}^{6}{B_{i}}$\ \par
with\ \par
$\displaystyle B_{1}=\frac{\displaystyle k^{2}}{\displaystyle 2}\displaystyle \int_{\Omega 
}^{}{\frac{\displaystyle \partial v}{\displaystyle \partial x}\frac{\displaystyle 
\partial }{\displaystyle \partial x}(x^{2}v)dxdy}$ ;  $\displaystyle 
B_{2}=\frac{\displaystyle 1}{\displaystyle 2}\displaystyle \int_{\Omega 
}^{}{x^{2}\frac{\displaystyle \partial v}{\displaystyle \partial y}\frac{\displaystyle 
\partial }{\displaystyle \partial y}(y^{2}v)dxdy}$;  $\displaystyle B_{3}=\rho k\displaystyle \int_{\Omega }^{}{x^{2}\frac{\displaystyle \partial 
v}{\displaystyle \partial x}}\ \frac{\displaystyle \partial }{\displaystyle 
\partial y}(yv)\ dxdy;$\ \par
\ \par
$\displaystyle B_{4}=-\displaystyle \int_{\Omega }^{}{\displaystyle \left({\displaystyle 
\alpha k^{2}x^{2}+g_{2}(x))+(1-\rho ^{2})k^{2}x^{2}\frac{\displaystyle 
\partial u_{1}}{\displaystyle \partial x}}\right) v\frac{\displaystyle 
\partial v}{\displaystyle \partial x}dxdy}$;$\displaystyle \ B_{5}=\ -\rho \alpha k\displaystyle \int_{\Omega }^{}{x^{2}y\frac{\displaystyle 
\partial v}{\displaystyle \partial y}vdxdy};$\ \par
\ \par
$\displaystyle B_{6}=-\displaystyle \int_{\Omega }^{}{\displaystyle \left({\displaystyle 
\frac{\displaystyle \alpha ^{2}k^{2}}{\displaystyle 2}x^{2}+\alpha g_{2}(x)+(1-\rho 
^{2})\alpha k^{2}x^{2}\frac{\displaystyle \partial u_{1}}{\displaystyle 
\partial x}}\right) v^{2}dxdy}.$\ \par
\ \par
It is easily seen, for $\epsilon _{i}>0,\ 1\leq i\leq 4$\ \par
\ \par
$\displaystyle B_{1}=\frac{\displaystyle k^{2}}{\displaystyle 2}\displaystyle \left\Vert{\displaystyle 
x\frac{\displaystyle \partial v}{\displaystyle \partial x}}\right\Vert 
^{2}_{\displaystyle L^{2}(\Omega )}-\frac{\displaystyle k^{2}}{\displaystyle 
2}\displaystyle \left\Vert{\displaystyle v}\right\Vert ^{2}_{\displaystyle 
L^{2}(\Omega )},$\ \par
\ \par
\ \par
$\displaystyle B_{2}=\frac{\displaystyle 1}{\displaystyle 2}\displaystyle \left\Vert{\displaystyle 
xy\frac{\displaystyle \partial v}{\displaystyle \partial y}}\right\Vert 
^{2}_{\displaystyle L^{2}(\Omega )}-\frac{\displaystyle 1}{\displaystyle 
2}\displaystyle \left\Vert{\displaystyle xv}\right\Vert ^{2}_{\displaystyle 
L^{2}(\Omega )},$\ \par
\ \par
\ \par
$\displaystyle B_{3}=\rho k\displaystyle \int_{\Omega }^{}{(x^{2}y\frac{\displaystyle 
\partial v}{\displaystyle \partial x}\frac{\displaystyle \partial v}{\displaystyle 
\partial y}-xv^{2})dxdy}$  and this term is bounded from below:\ \par
\ \par
$\displaystyle \displaystyle \left\vert{\displaystyle B_{3}}\right\vert \geq -\frac{\displaystyle 
\displaystyle \left\vert{\displaystyle \rho }\right\vert }{\displaystyle 
2}\displaystyle \left({\displaystyle \displaystyle \left\Vert{\displaystyle 
xy\frac{\displaystyle \partial v}{\displaystyle \partial y}}\right\Vert 
^{2}_{\displaystyle L^{2}(\Omega )}+k^{2}\displaystyle \left\Vert{\displaystyle 
x\frac{\displaystyle \partial v}{\displaystyle \partial x}}\right\Vert 
^{2}_{\displaystyle L^{2}(\Omega )}}\right) -\displaystyle \left\vert{\displaystyle 
\rho }\right\vert \displaystyle \left({\displaystyle \epsilon _{1}\displaystyle 
\left\Vert{\displaystyle xv}\right\Vert ^{2}_{\displaystyle L^{2}(\Omega 
)}+\frac{\displaystyle k^{2}}{\displaystyle 4\epsilon _{1}}\displaystyle 
\left\Vert{\displaystyle v}\right\Vert ^{2}_{\displaystyle L^{2}(\Omega 
)}}\right) $, $\epsilon _{1}>0,$\ \par
\ \par
\ \par
$\displaystyle B_{4}=\frac{\displaystyle 1}{\displaystyle 2}\displaystyle \int_{\Omega 
}^{}{v^{2}\displaystyle \left({\displaystyle 2\alpha k^{2}x+g_{2}'(x)}\right) 
dx}+(1-\rho ^{2})k^{2}\displaystyle \int_{\Omega }^{}{\frac{\displaystyle 
\partial u_{1}}{\displaystyle \partial x}x^{2}v\frac{\displaystyle \partial 
v}{\displaystyle \partial x}dxdy}.$\ \par
\ \par
This term is bounded from below:\ \par
\ \par
$\displaystyle B_{4}\geq \displaystyle \int_{\Omega }^{}{\displaystyle \left({\displaystyle 
\alpha k^{2}-\delta }\right) xv^{2}dxdy}-\frac{\displaystyle 1}{\displaystyle 
2}\displaystyle \left\vert{\displaystyle \rho }\right\vert k\displaystyle 
\left\Vert{\displaystyle f'_{1}}\right\Vert _{\displaystyle L^{\infty 
}({\mathbb{R}}_{+})}\displaystyle \left\Vert{\displaystyle v}\right\Vert 
^{2}_{\displaystyle L^{2}(\Omega )}$\ \par
$\displaystyle -(1-\rho ^{2})\displaystyle \left({\displaystyle \epsilon _{3}\displaystyle 
\left\Vert{\displaystyle x\frac{\displaystyle \partial v}{\displaystyle 
\partial x}}\right\Vert ^{2}_{\displaystyle L^{2}(\Omega )}+\frac{\displaystyle 
k^{4}}{\displaystyle 4\epsilon _{3}}\displaystyle \left\Vert{\displaystyle 
x\frac{\displaystyle \partial u_{1}}{\displaystyle \partial x}}\right\Vert 
^{2}_{\displaystyle L^{\infty }(\Omega )}\displaystyle \left\Vert{\displaystyle 
v}\right\Vert ^{2}_{\displaystyle L^{2}(\Omega )}}\right) $\ \par
\ \par
and we obtain\ \par
$\displaystyle B_{4}\geq -\epsilon _{2}\displaystyle \left\Vert{\displaystyle xv}\right\Vert 
^{2}_{\displaystyle L^{2}(\Omega )}-(1-\rho ^{2})\epsilon _{3}\displaystyle 
\left\Vert{\displaystyle x\frac{\displaystyle \partial v}{\displaystyle 
\partial x}}\right\Vert ^{2}_{\displaystyle L^{2}(\Omega )}$\ \par
$\displaystyle -\displaystyle \left({\displaystyle \frac{\displaystyle (\alpha k^{2}-\delta 
)^{2}}{\displaystyle 4\epsilon _{2}}-\frac{\displaystyle 1}{\displaystyle 
2}\displaystyle \left\vert{\displaystyle \rho }\right\vert k\displaystyle 
\left\Vert{\displaystyle f'_{1}}\right\Vert _{\displaystyle L^{\infty 
}({\mathbb{R}}_{+})}+\frac{\displaystyle k^{4}}{\displaystyle 4\epsilon 
_{3}}\displaystyle \left\Vert{\displaystyle x\frac{\displaystyle \partial 
u_{1}}{\displaystyle \partial x}}\right\Vert ^{2}_{\displaystyle L^{\infty 
}(\Omega )}}\right) \displaystyle \left\Vert{\displaystyle v}\right\Vert 
^{2}_{\displaystyle L^{2}(\Omega )}.$\ \par
\ \par
\ \par
$\displaystyle B_{5}\geq -\frac{\displaystyle \displaystyle \left\vert{\displaystyle 
\rho }\right\vert }{\displaystyle 2}\alpha k\displaystyle \left\Vert{\displaystyle 
xv}\right\Vert ^{2}_{\displaystyle L^{2}(\Omega )}$.\ \par
\ \par
\ \par
$\displaystyle B_{6}\geq \alpha \displaystyle \left({\displaystyle \delta -\frac{\displaystyle 
\alpha k^{2}}{\displaystyle 2}}\right) \displaystyle \left\Vert{\displaystyle 
xv}\right\Vert ^{2}_{\displaystyle L^{2}(\Omega )}-\displaystyle \int_{\Omega 
}^{}{\alpha (\delta \sigma _{1}+(1-\rho ^{2})k^{2}\displaystyle \left\Vert{\displaystyle 
x\frac{\displaystyle \partial u_{1}}{\displaystyle \partial x}}\right\Vert 
^{2}_{\displaystyle L^{\infty }(\Omega )})xv^{2}dxdy}$\ \par
$\displaystyle -\alpha \displaystyle \left\vert{\displaystyle \rho }\right\vert k\displaystyle 
\left\Vert{\displaystyle f_{1}}\right\Vert _{\displaystyle L^{\infty 
}({\mathbb{R}}_{+})}\displaystyle \left\Vert{\displaystyle v}\right\Vert 
^{2}_{\displaystyle L^{2}(\Omega )}.$\ \par
\ \par
and we deduce :\ \par
$\displaystyle B_{6}\geq \displaystyle \left({\displaystyle \alpha \displaystyle \left({\displaystyle 
\delta -\frac{\displaystyle \alpha k^{2}}{\displaystyle 2}}\right) -\epsilon 
_{4}}\right) \displaystyle \left\Vert{\displaystyle xv}\right\Vert ^{2}$\ 
\par
$\displaystyle -\displaystyle \left({\displaystyle \frac{\displaystyle \alpha ^{2}}{\displaystyle 
4\epsilon _{4}}(\delta \sigma _{1}+(1-\rho ^{2})k^{2}\displaystyle \left\Vert{\displaystyle 
x\frac{\displaystyle \partial u_{1}}{\displaystyle \partial x}}\right\Vert 
^{2}_{\displaystyle L^{\infty }(\Omega )})^{2}+\alpha \displaystyle \left\vert{\displaystyle 
\rho }\right\vert k\displaystyle \left\Vert{\displaystyle f_{1}}\right\Vert 
_{\displaystyle L^{\infty }({\mathbb{R}}_{+})}}\right) \displaystyle 
\left\Vert{\displaystyle v}\right\Vert ^{2}_{\displaystyle L^{2}(\Omega 
)}.$\ \par
\ \par
\ \par
From all these estimates, we get:\ \par
$\displaystyle \hat b(v,v)\geq (1-\displaystyle \left\vert{\displaystyle \rho }\right\vert 
)\frac{\displaystyle k^{2}}{\displaystyle 2}(1-(1+\displaystyle \left\vert{\displaystyle 
\rho }\right\vert )\epsilon _{3})\displaystyle \left\Vert{\displaystyle 
x\frac{\displaystyle \partial v}{\displaystyle \partial x}}\right\Vert 
^{2}_{\displaystyle L^{2}(\Omega )}+\frac{\displaystyle 1}{\displaystyle 
2}(1-\displaystyle \left\vert{\displaystyle \rho }\right\vert )\displaystyle 
\left\Vert{\displaystyle xy\frac{\displaystyle \partial v}{\displaystyle 
\partial y}}\right\Vert ^{2}_{\displaystyle L^{2}(\Omega )}$\ \par
$\displaystyle +\displaystyle \left({\displaystyle \alpha \displaystyle \left({\displaystyle 
\delta -\frac{\displaystyle \alpha k^{2}}{\displaystyle 2}}\right) -\frac{\displaystyle 
1}{\displaystyle 2}-\displaystyle \left\vert{\displaystyle \rho }\right\vert 
\frac{\displaystyle \alpha k}{\displaystyle 2}-\epsilon }\right) \displaystyle 
\left\Vert{\displaystyle xv}\right\Vert ^{2}_{\displaystyle L^{2}(\Omega 
)}$ $-c\left\Vert{v}\right\Vert ^{2}_{L^{2}(\Omega )}$\ \par
\ \par
with $\epsilon =\epsilon _{1}+\epsilon _{2}+\epsilon _{3}+\epsilon _{4}$ 
and $c$ is a positive constant depending on the different parameters.\ 
\par
\ \par
Let us denote  $\displaystyle P(\alpha )=\alpha \displaystyle \left({\displaystyle \delta -\frac{\displaystyle 
\alpha k^{2}}{\displaystyle 2}}\right) -\frac{\displaystyle 1}{\displaystyle 
2}-\displaystyle \left\vert{\displaystyle \rho }\right\vert \frac{\displaystyle 
\alpha k}{\displaystyle 2}-\epsilon .$\ \par
This polynomial admits two positive real roots if $k<\frac{2}{3}\delta $ 
and $\epsilon $ small enough. This last condition is generally satisfied 
in practice ($k=0.4,\ \delta =2)$. So we can choose $\alpha >0$ such 
that $P(\alpha )>0$.\ \par
Therefore by choosing $\epsilon _{i},1\leq i\leq 4$ small enough, $\alpha 
$ such that $P(\alpha )>0$ and if $\left\vert{\rho }\right\vert <1$, 
we obtain: $\displaystyle \hat b(v,v)\geq C\displaystyle \left\Vert{\displaystyle v}\right\Vert 
_{2}^{2}-c\displaystyle \left\Vert{\displaystyle v}\right\Vert ^{2}_{\displaystyle 
L^{2}(\Omega )}$.\ \par
\ \par
We may write problem ~(\ref{1MarchesIncomplets0}) in variational form:\ 
\par

\begin{equation} 
\displaystyle \left\lbrace{\displaystyle \begin{matrix}{\displaystyle 
Find\ \hat u_{2}\in C(0,T;L^{2}(\Omega ))\cap L^{2}(0,T;\hat V_{2})}\cr 
{\displaystyle (\frac{\displaystyle \partial \hat u_{2}}{\displaystyle 
\partial t},v)+\hat b(\hat u_{2},v)=0,\ \forall v\in V_{2}}\cr {\displaystyle 
\hat u_{2}(0)=e^{-\alpha x}F}\cr \end{matrix}}\right. .\label{1MarchesIncomplets2}
\end{equation} \ \par
\ \par
Since $\hat u_{2}(0)\in \hat V_{2}$, by using proposition ~\ref{1MarchesIncomplets1}, 
we deduce the following theorem ~\cite{Lions}:\ \par
\ \par
\begin{Theorem} Under the hypotheses of Proposition ~\ref{1MarchesIncomplets1}, 
problem ~(\ref{1MarchesIncomplets2}) admits a unique solution\ \par
\end{Theorem}
\ \par
\subsection{Numerical solution}
\ \par
For the computation of $u_{2},$, we shall use a backward Euler method 
in time and a finite element method in space~\cite{Ciarlet}.\ \par
\ \par
We define a triangulation ${\mathcal{T}}_{h}$ of $\Omega $ in the following 
way:\ \par
- $(x_{i}),\ 0\leq i\leq N$ is the sequence defined in section (2.1.2).\ 
\par
- Let $(y_{j}),\ 0\leq j\leq M$ ( $y_{0}=0$) another increasing sequence.\ 
\par
We denote $k_{j}=y_{j}-y_{j-1}$ and we assume that the sequence $(k_{j}),\ 1\leq j\leq M$ 
is increasing.\ \par
\ \par
The domain $[0,x_{N})\times [0,y_{M}]$ is divided in rectangles $(x_{i-1},x_{i})\times (y_{j-1},y_{j})$, 
$1\leq i\leq N,\ 1\leq j\leq M$ and each rectangle is divided in two 
triangles by the first diagonal. We denote ${\mathcal{T}}_{h}^{1}$ the 
set of these triangles.\ \par
We define ${\mathcal{T}}_{h}^{2},\ {\mathcal{T}}_{h}^{3},\ K_{NM}$ by:\ 
\par
$\displaystyle {\mathcal{T}}_{h}^{2}=\displaystyle \left\lbrace{\displaystyle \displaystyle 
\left.{\displaystyle K_{j}/K_{j}=(x_{N},+\infty [\times (y_{j-1},y_{j}),\ 
j=1,M}\right\rbrace }\right. $\ \par
\ \par
$\displaystyle {\mathcal{T}}_{h}^{3}=\displaystyle \left\lbrace{\displaystyle \displaystyle 
\left.{\displaystyle K_{i}/K_{i}=(x_{i-1},x_{i})\times (y_{M},+\infty 
),i=1,N}\right\rbrace }\right. $\ \par
\ \par
$\displaystyle K_{NM}=(x_{N},+\infty )\times (y_{M},+\infty )$\ \par
\ \par
The triangulation is then defined by: $\displaystyle {\mathcal{T}}_{h}={\mathcal{T}}_{h}^{1}\cup {\mathcal{T}}_{h}^{2}\cup 
{\mathcal{T}}_{h}^{3}\cup K_{NM}$.\ \par
\ \par
We associate to this triangulation the finite-dimensional space $V_{2h}\ $defined 
by:\ \par
\ \par
$\displaystyle V_{2h}=\displaystyle \left\lbrace{\displaystyle \displaystyle \left.{\displaystyle 
\begin{matrix}{\displaystyle v_{h}\in C^{0}(\Omega )/\forall K\in {\mathcal{T}}_{h}^{1},\ 
v_{\displaystyle h\mid K}\in P_{1},\ \forall K\in {\mathcal{T}}_{h}^{2},\ 
v_{\displaystyle h\mid K}\in P_{y1}}\cr {\displaystyle \forall K\in {\mathcal{T}}_{h}^{3},\ 
yv_{\displaystyle h\mid K}\in P_{x1},\ yv_{\displaystyle h\mid K_{NM}}\in 
P_{0}}\cr \end{matrix}}\right\rbrace }\right. $\ \par
\ \par
\ \par
$P_{1}$ is the space of polynomials of degree $\leq 1$ in $x,y$ ;$\ P_{y1}$ 
is the space of polynomials of degree $\leq 1$ in $y$ ;$\ P_{x1}$ is 
the space of polynomials of degree $\leq 1$ in $x$ ;$\ P_{0}$ is the 
space of constants.\ \par
\ \par
Let $\alpha >0$ such that $P(\alpha )>0$.\ \par
If $v_{h}\in V_{2h}$, it is easily seen that $\hat v_{h}=e^{-\alpha x}v_{h}\in \hat V_{2}$.\ 
\par
The approximate value of $u_{2}$ at the time level $t_{n}$  will be in$\ V_{2h}$.\ 
\par
\ \par
If $v_{h}\in V_{2h}$, we denote $\tilde v_{h}=e^{-2\alpha x}v_{h}$ and 
we define the bilinear form $b$ on $V_{2h}\times V_{2h}$ by: $\forall u_{h},v_{h}\in V_{2h},$ 
\ \par
$\displaystyle b(u_{h},v_{h})=\frac{\displaystyle k^{2}}{\displaystyle 2}\displaystyle 
\int_{\Omega }^{}{\frac{\displaystyle \partial u_{h}}{\displaystyle \partial 
x}\frac{\displaystyle \partial }{\displaystyle \partial x}(x^{2}\tilde 
v_{h})dxdy}+\frac{\displaystyle 1}{\displaystyle 2}\displaystyle \int_{\Omega 
}^{}{x^{2}\frac{\displaystyle \partial u_{h}}{\displaystyle \partial 
y}\frac{\displaystyle \partial }{\displaystyle \partial y}(y^{2}\tilde 
v_{h})dxdy}$\ \par
$\displaystyle +\rho k\displaystyle \int_{\Omega }^{}{x^{2}\frac{\displaystyle \partial 
u_{h}}{\displaystyle \partial x}\frac{\displaystyle \partial }{\displaystyle 
\partial y}(y\tilde v_{h})dxdy}$ $\displaystyle -\displaystyle \int_{\Omega }^{}{xg_{2}(x)\frac{\displaystyle \partial 
u_{h}}{\displaystyle \partial x}\tilde v_{h}dxdy}$\ \par
$\displaystyle -(1-\rho ^{2})k^{2}\displaystyle \int_{\Omega }^{}{x^{2}\frac{\displaystyle 
\partial u_{1}}{\displaystyle \partial x}\frac{\displaystyle \partial 
u_{h}}{\displaystyle \partial x}\tilde v_{h}dxdy}$.\ \par
\ \par
We have the equality: $b(u_{h},v_{h})=\hat b(\hat u_{h},\hat v_{h}),\ \forall u_{h},v_{h}\in 
V_{2h}$  and we define the approximate value $u_{2h}^{n+1}$ of $u_{2}(t_{n+1})$ 
as the solution of the following problem:\ \par
\ \par
$\displaystyle (u_{2h}^{n+1}-u_{2h}^{n},\tilde v_{h})_{h}+b_{h}(u_{2h}^{n+1},v_{h})=0,$ 
$\forall v_{h}\in V_{2h}$,\ \par
$u_{2h}^{0}=F_{h},$\ \par
\ \par
where $F_{h}$ is the Lagrange interpolate of $F$ in $V_{2h}$; \ \par
$(u_{h},v_{h})_{h}$   is an approximate scalar product in $L^{2}(\Omega )$ 
and $b_{h}$ an approximation of $b$, obtained by using numerical integration.\ 
\par
\ \par
In Fig 2, we present the variation of $u_{2}(T)$  in $y\ $(or $b(0)\ $in 
$P$) for different values of the volatility; the parameter $\rho \ $ 
is null; the exercise price $E$ is equal to $1$; the others parameters 
have the same values as in Fig1.\ \par
\ \par

\begin{figure}[h]
\begin{center}
\rotatebox{0}{\resizebox{10cm}{!}{\includegraphics{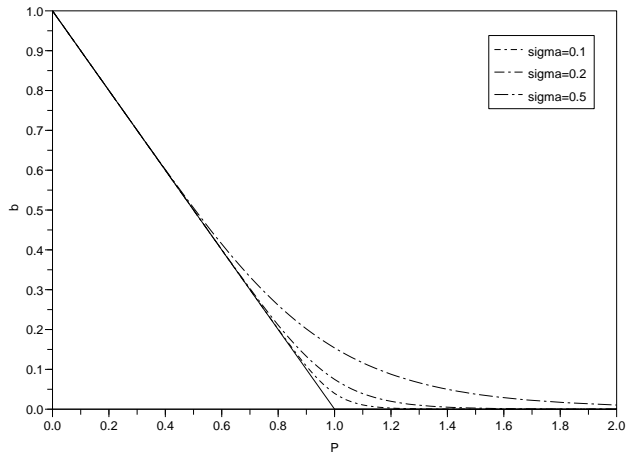}}}
\caption{b(0) with T=1}
\label{1MarchesIncomplets100}
\end{center}
\end{figure}
\ \par
\ \par
Fig 3 represents the variation of $u_{2}(T)\ \ $in $y\ $for different 
values of $\rho $, the volatilty is equal to 0.3. \ \par
\ \par

\begin{figure}[h]
\begin{center}
\rotatebox{0}{\resizebox{10cm}{!}{\includegraphics{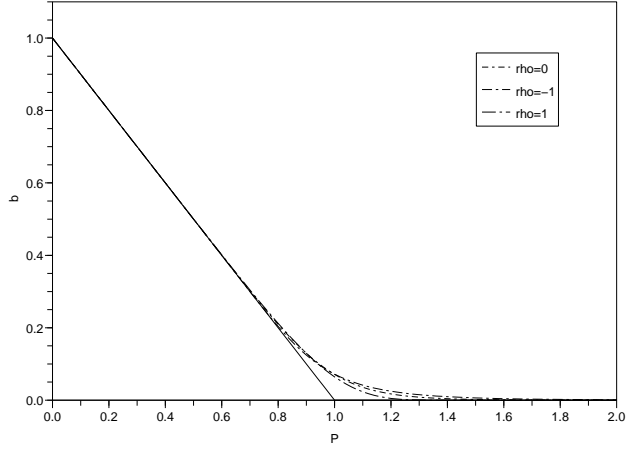}}}
\caption{b(0) with sigma=0.3}
\label{1MarchesIncomplets101}
\end{center}
\end{figure}
\ \par
\ \par
\ \par
\section{Computation of $u_{3}$}
\setcounter{equation}{0}$u_{3}$ is solution of \ \par
\ \par
$\displaystyle \frac{\displaystyle \partial u_{3}}{\displaystyle \partial t}-\frac{\displaystyle 
k^{2}x^{2}}{\displaystyle 2}\frac{\displaystyle \partial ^{2}u_{3}}{\displaystyle 
\partial x^{2}}-\frac{\displaystyle x^{2}y^{2}}{\displaystyle 2}\frac{\displaystyle 
\partial ^{2}u_{3}}{\displaystyle \partial y^{2}}-\rho kx^{2}y\frac{\displaystyle 
\partial ^{2}u_{3}}{\displaystyle \partial x\partial y}-g(x)\frac{\displaystyle 
\partial u_{3}}{\displaystyle \partial x}-xf(x)y\frac{\displaystyle \partial 
u_{3}}{\displaystyle \partial y}$\ \par
$\displaystyle =(1-\rho ^{2})k^{2}x^{2}exp(u_{1})\displaystyle \left({\displaystyle \frac{\displaystyle 
\partial u_{2}}{\displaystyle \partial x}}\right) ^{2},$\ \par
with the initial condition: $u_{3}(0)=0$.\ \par
\ \par
\subsection{Existence and uniqueness of the solution}
As for the computation of $u_{2}$, we make the change of unknown: $\hat u_{3}=e^{-\alpha x}u_{3}$ 
and $\hat u_{3}$ is solution of:\ \par
\ \par
$\displaystyle \frac{\displaystyle \partial \hat u_{3}}{\displaystyle \partial t}-\frac{\displaystyle 
k^{2}x^{2}}{\displaystyle 2}\frac{\displaystyle \partial ^{2}\hat u_{3}}{\displaystyle 
\partial x^{2}}-\frac{\displaystyle x^{2}y^{2}}{\displaystyle 2}\frac{\displaystyle 
\partial ^{2}\hat u_{3}}{\displaystyle \partial y^{2}}-\rho kx^{2}y\frac{\displaystyle 
\partial ^{2}\hat u_{3}}{\displaystyle \partial x\partial y}-\frac{\displaystyle 
\partial \hat u_{3}}{\displaystyle \partial x}\displaystyle \left({\displaystyle 
\alpha k^{2}x^{2}+g(x)}\right) -\frac{\displaystyle \partial \hat u_{3}}{\displaystyle 
\partial y}\displaystyle \left({\displaystyle \rho k\alpha x^{2}y+xf(x)y}\right) 
$\ \par

\begin{equation} 
-\hat u_{3}\displaystyle \left({\displaystyle \frac{\displaystyle \alpha 
^{2}k^{2}x^{2}}{\displaystyle 2}+\alpha g(x)}\right) =(1-\rho ^{2})k^{2}x^{2}exp(u_{1}-\alpha 
x)\displaystyle \left({\displaystyle \frac{\displaystyle \partial u_{2}}{\displaystyle 
\partial x}}\right) ^{2}.\label{1MarchesIncomplets3}
\end{equation} \ \par
\ \par
We define the bilinear form $c$ on $\hat V_{2}\times \hat V_{2}$ by:\ 
\par
$\displaystyle \forall u,v\in \hat V_{2},\ \hat c(u,v)=\frac{\displaystyle k^{2}}{\displaystyle 
2}\displaystyle \int_{\Omega }^{}{\frac{\displaystyle \partial u}{\displaystyle 
\partial x}\frac{\displaystyle \partial }{\displaystyle \partial x}(x^{2}v)dxdy}+\frac{\displaystyle 
1}{\displaystyle 2}\displaystyle \int_{\Omega }^{}{x^{2}\frac{\displaystyle 
\partial u}{\displaystyle \partial y}\frac{\displaystyle \partial }{\displaystyle 
\partial y}(y^{2}v)dxdy}$\ \par
$\displaystyle +\rho k\displaystyle \int_{\Omega }^{}{x^{2}\frac{\displaystyle \partial 
u}{\displaystyle \partial x}\frac{\displaystyle \partial }{\displaystyle 
\partial y}(yv)dxdy}-\displaystyle \int_{\Omega }^{}{(\alpha k^{2}x^{2}+g(x))\frac{\displaystyle 
\partial u}{\displaystyle \partial x}vdxdy}$\ \par
$\displaystyle -\displaystyle \int_{\Omega }^{}{(\rho \alpha kx^{2}+xf(x))y\frac{\displaystyle 
\partial u}{\displaystyle \partial y}vdxdy}$ $\displaystyle -\displaystyle \int_{\Omega }^{}{(\frac{\displaystyle \alpha ^{2}k^{2}x^{2}}{\displaystyle 
2}+\alpha g(x))uvdxdy}.$\ \par
\ \par
\ \par
It is clear that this bilinear form is continue on $\hat V_{2}\times \hat V_{2}$ 
and under the same hypotheses  as for $b$, it satisfies: $\displaystyle 
\ \forall v\in \hat V_{2},\ \hat c(v,v)\geq C\displaystyle \left\Vert{\displaystyle 
v}\right\Vert ^{2}_{\displaystyle \hat V_{2}}-c\displaystyle \left\Vert{\displaystyle 
v}\right\Vert ^{2}_{\displaystyle L^{2}(\Omega )}.\ $\ \par
\ \par
We denote  $\displaystyle G(x)=(1-\rho ^{2})k^{2}x^{2}e^{\displaystyle u_{1}}\displaystyle \left({\displaystyle 
\frac{\displaystyle \partial u_{2}}{\displaystyle \partial x}}\right) 
^{2},\ x>0$ . \ \par
Problem ~(\ref{1MarchesIncomplets3}) may be written in variational form:\ 
\par

\begin{equation} 
\displaystyle \left\lbrace{\displaystyle \begin{matrix}{\displaystyle 
Find\ \ \hat u_{3}\in C^{2}(0,T;L^{2}(\Omega )\cap L^{2}(0,T;\hat V_{2})}\cr 
{\displaystyle \displaystyle \left({\displaystyle \frac{\displaystyle 
\partial \hat u_{3}}{\displaystyle \partial t},v}\right) +\hat c(\hat 
u_{3},v)=(Ge^{-\alpha x},v)}\cr {\displaystyle \hat u_{3}(0)=0}\cr \end{matrix}}\right. 
{\mathrm{.}}\label{1MarchesIncomplets4}
\end{equation} \ \par
\ \par
If we assume that $\frac{\partial u_{2}}{\partial x}\in C(0,T;L^{\infty }(\Omega ))$, 
the second member is continue on $\hat V_{2}$ and we get the following 
theorem:\ \par
\ \par
\begin{Theorem} If  $\left\vert{\rho }\right\vert <1$, $k<\frac{2}{3}\delta $ 
and $\alpha $ such that $\hat c$ satisfies Garding inequality, problem 
~(\ref{1MarchesIncomplets4}) admits a unique solution.\ \par
\end{Theorem}
\ \par
\subsection{Numerical solution}
To compute $u_{3}$, we use the same method as for the computation of 
$u_{2}$.\ \par
\ \par
If $v_{h}\in V_{2h}$, we denote $\displaystyle \tilde v_{h}=e^{-2\alpha x}v_{h}$ 
and define the bilinear form on $V_{2h}\times V_{2h}$ by:\ \par
\ \par
$\displaystyle c(u_{h},v_{h})=\frac{\displaystyle k^{2}}{\displaystyle 2}\displaystyle 
\int_{\Omega }^{}{\frac{\displaystyle \partial u_{h}}{\displaystyle \partial 
x}\frac{\displaystyle \partial (x^{2}\tilde v_{h})}{\displaystyle \partial 
x}dxdy}+\frac{\displaystyle 1}{\displaystyle 2}\displaystyle \int_{\Omega 
}^{}{x^{2}\frac{\displaystyle \partial u_{h}}{\displaystyle \partial 
y}\frac{\displaystyle \partial (y^{2}\tilde v_{h})}{\displaystyle \partial 
y}dxdy}+\rho k\displaystyle \int_{\Omega }^{}{x^{2}\frac{\displaystyle 
\partial u_{h}}{\displaystyle \partial x}\frac{\displaystyle \partial 
(y\tilde v_{h})}{\displaystyle \partial y}dxdy}$\ \par
$\displaystyle -\displaystyle \int_{\Omega }^{}{g(x)\frac{\displaystyle \partial u_{h}}{\displaystyle 
\partial x}\tilde v_{h}dxdy}-\displaystyle \int_{\Omega }^{}{xf(x)y\frac{\displaystyle 
\partial u_{h}}{\displaystyle \partial x}\tilde v_{h}dxdy}$ \ \par
\ \par
and we have the equality: $\displaystyle \forall u_{h},v_{h}\in V_{2h},\ c(u_{h},v_{h})=\hat c(\hat u_{h},\hat 
v_{h})$.\ \par
\ \par
The approximate solution $u_{3h}^{n+1}\ $ at the time level $t_{n+1}$ 
satisfies:\ \par
\ \par
$\displaystyle (u_{3h}^{n+1}-u_{3h}^{n},\tilde v_{h})_{h}+\Delta t_{n}c_{h}(u_{3h}^{n+1},vh)=(G,\tilde 
v_{h})_{h},\ \forall v_{h}\in V_{2h},$\ \par
$\displaystyle u_{3h}^{0}=0.$\ \par
\ \par
where $c_{h}$ is an approximation of $c$ obtained by using numerical 
integration.\ \par
\ \par
Fig4 represents the variation of$\ u_{3}(T)$ in $y\ $(or $c(0)$ in $P$) 
for different values of $\sigma $. The parameter $\rho $ is equal to 
$0$. (In the case $\rho =1$ or $\rho =-1$, $u_{3}$ is null). The replication 
error $\epsilon ^{*}$ is given by: $\epsilon ^{*}=\sqrt{c(0)}$.\ \par
\ \par

\begin{figure}[h]
\begin{center}
\rotatebox{0}{\resizebox{10cm}{!}{\includegraphics{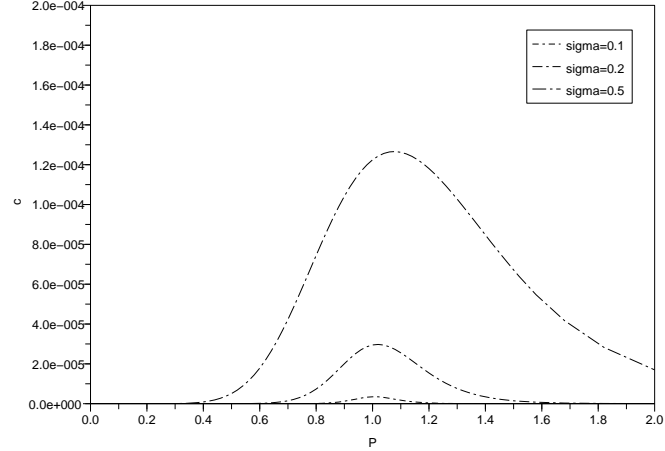}}}
\caption{c(0) with T=1}
\label{1MarchesIncomplets102}
\end{center}
\end{figure}
\ \par
\ \par
\ \par
Knowing $u_{2}$, we can compute the least-cost optimal strategy $\displaystyle 
\theta ^{*}(0)=\frac{\displaystyle \partial b}{\displaystyle \partial 
P}(0)+\frac{\displaystyle \rho k}{\displaystyle P}\frac{\displaystyle 
\partial b}{\displaystyle \partial \sigma }(0)=\frac{\displaystyle \partial 
u_{2}}{\displaystyle \partial y}(T)+\frac{\displaystyle \rho k}{\displaystyle 
y}\frac{\displaystyle \partial u_{2}}{\displaystyle \partial y}(T)$.\ 
\par
\ \par
Fig 5 represents the variation of $\theta ^{*}$ in $P$ for different 
values of $\sigma $; the parameter $\rho $ is equal to 0.\ \par

\begin{figure}[h]
\begin{center}
\rotatebox{0}{\resizebox{10cm}{!}{\includegraphics{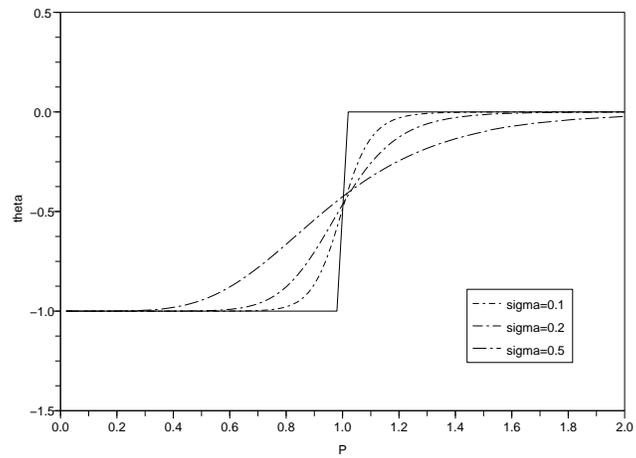}}}
\caption{Optimal strategy}
\label{1MarchesIncomplets103}
\end{center}
\end{figure}
\ \par

 \ \par
\ \par
\ \par

\bibliographystyle{C:/TexLive/texmf/bibtex/bst/base/plain}

\begin{thebibliography}{1}

\bibitem{Kogan}
D.~Bertsimas, L.~Kogan, and A.W. Lo.
\newblock {Hedging derivative securities and incomplete markets: An
  $\epsilon$-arbitrage approach}.
\newblock {\em Operations Research}, 49(3):372--397, 2001.

\bibitem{Black}
F.~Black and M.~Scholes.
\newblock Pricing of options and corporate liabilities.
\newblock {\em J. of Political Econom.}, 81:637--654, 1973.

\bibitem{Bolley}
P.~Bolley and J.~Camus J.
\newblock {Quelques r\'esultats sur les espaces de Sobolev avec poids}.
\newblock {\em {Publications de S\'eminaires Math\'ematiques de Rennes}},
  1:1--69, 1969.

\bibitem{Ciarlet}
P.G. Ciarlet.
\newblock {\em The Finite Element Method for Elliptic Problems}.
\newblock North-Holland Publishing Company, 1978.

\bibitem{Lions}
J.L. Lions and E.~Magenes.
\newblock {\em Problemes aux limites non homogènes et applications}, volume~1.
\newblock Dunod, Paris, 1968.

\bibitem{Ritchmyer}
R.D. Ritchmyer and K.W. Morton.
\newblock {\em Difference Methods for Initial-Value Problems}.
\newblock 2nd ed, Wiley-Interscience, New York, 1967.

\bibitem{Simon}
J.~Simon.
\newblock {Compact sets in the space $L^{p}(0,T;B)$}.
\newblock {\em Ann.Mat.Pura.Appl}, 146:65--96, 1987.

\end{thebibliography}

\end{document}